\documentclass[12 pt]{amsart}

\usepackage{amscd,amsfonts,amssymb,amsmath}
\usepackage{tikz-cd}
\usepackage{xcolor}
\usepackage{placeins}
\usepackage{hyperref}

\newtheorem{theorem}{Theorem}[section]
\theoremstyle{definition}

\newtheorem{lemma}[theorem]{Lemma}
\newtheorem{proposition}[theorem]{Proposition}
\newtheorem{question}[theorem]{Question}

\newtheorem{definition}[theorem]{Definition}
\newtheorem{example}[theorem]{Example}

\numberwithin{equation}{section}

\newcommand{\id}{\mathrm{id}}

\title{Lie  Rota--Baxter operators on the Sweedler algebra~$H_4$
}
\author[Valeriy G.~Bardakov,  Igor M.~Nikonov, and Viktor N. Zhelaybin]{Valeriy G.~Bardakov,  Igor M.~Nikonov, and Viktor N. Zhelaybin}

\date{\today}

\begin{document}
\sloppy



\maketitle
\begin{abstract}
If $A$ is an associative algebra, then we can define the adjoint  Lie algebra $A^{(-)}$ and Jordan  algebra $A^{(+)}$.   It is easy to see that  any associative Rota--Baxter operator on $A$ induces  a Lie and Jordan  Rota--Baxter operator on $A^{(-)}$ and $A^{(+)}$ respectively. Are there Lie  (Jordan)   Rota--Baxter operators, which are not associative Rota--Baxter operators?

In the present article we are studying these questions for the Sweedler algebra $H_4$, that is a 4-dimension non-commutative Hopf algebra.  More precisely, we describe the Rota--Baxter operators on Lie algebra on the adjoint Lie algebra $H_4^{(-)}$.


\end{abstract}

\maketitle
\tableofcontents

\section{Introduction}

Rota--Baxter operators on algebras are known since the middle of the previous century 
and they have in turn connections with mathematical physics,
number theory, operad theory, Hopf algebras, see the monograph~\cite{GuoMonograph}.
The Rota--Baxter  operators (RB-operator, for short) first appeared in the work \cite{Tricomi}. These
operators became popular after the work~\cite{Baxter},  where they
were  defined as formal generalization of integration by parts formula.
  Further,  these operators  were  studied  on commutative algebras in \cite{Atkinson, Cartier, Rota}.

At the beginning of the 1980 years the notion to the classical
Yang--Baxter equation was introduced \cite{Drinfeld}. The solutions of this equation have a deep connection  with RB-operators \cite{STyanSh}. Actually, in \cite{STyanSh} the notion of an RB-operator of non zero weight on a Lie algebra was introduced. Such operators  are  called  solutions to a modified Yang--Baxter equation.

One of the interesting direction in the study of Rota--Baxter operators is a problem
of classification of RB-operators on a given algebra. The  RB-operators were classified on
$sl_2(\mathbb{C})$ in \cite{Kolesniukov, Konovalova, Pan, PeiBai}, on $M_2(\mathbb{C})$ in \cite{BG, Rota}, on $sl_3(\mathbb{C})$  in \cite{Konovalova} and on other algebras  in \cite{AnBai, BGP, DuBaiGou}.
The  skew-symmetric RB-operators of nonzero weight on
$M_3(\mathbb{C})$ were  described in \cite{Sokolov}.
Up to conjugation with automorphisms and transpose, 8 series of RB-operators were obtained in \cite{Sokolov}.
All RB-operators of nonzero weight on the matrix algebra
of order three over an algebraically closed field of characteristic zero were classified in \cite{GonGub}. In \cite{BG}, the Rota–Baxter operators of weight zero or nonzero on the simple four-dimensional Jordan  superalgebra $D_t$ over an algebraically closed field of characteristic~0 was described.
 A classification of RB-operators of
weight 1 on the general linear complex Lie algebra of order 2 was given in \cite{Goncharov,GK}.  The RB-operators on four-dimensional Sweedler's algebra considered as an associative algebra were founded in \cite{Ma-1}.

In the work \cite{GuoLShen}, the notion of a Rota--Baxter operator for groups was introduced.  A~mapping $B \colon  G \to  G$ of a group $G$ is called a Rota--Baxter operator on $G$
if
$$B(g)B(h) = B(gB(g)hB(g)^{-1})$$
 for any $g, h \in  G$.
In the same
paper, it was proved that if $(G, B)$ is a Rota--Baxter Lie group, then the tangent map
of $B$ at the identity is a RB-operator of weight 1 on the Lie algebra of the Lie group $G$. One can find different constructions of RB-operators on groups in \cite{BG-2}. In \cite{BG-1}, a connection between groups with RB-operators and skew braces was found. Important progress in this direction has been made in the work \cite{Goncharov2021}, where a concept of a Rota--Baxter operator was introduced and studied on a cocommutative Hopf algebra. It should also be noted that the notion of Rota--Baxter operator on a cocommutative Hopf algebra generalizes the notion of a Rota--Baxter operator on a group.

Let $H$ be a cocommutative Hopf algebra. Then a Rota--Baxter operator on the Hopf algebra $H$ is a mapping
$B \colon H \to H$, which is a coalgebra homomorphism
and satisfies the identity
\begin{equation}\label{0} B(a)B(b) = B(a_{(1)}B(a_{(2)}bS(a_{(3)})) \, a,b\in H,\end{equation} where $\Delta(a)=a_{(1)}\otimes a_{(2)}$  and $a_{(1)}\otimes a_{(2)}\otimes a_{(3)}=(\Delta\otimes id)\Delta(a)$ are Sweedler's notation. A definition of RB-operators on arbitrary Hopf algebras can be found in \cite{BN}. In \cite{ZLMZ}, Hopf braces were defined and studied.

It was shown in \cite{BNZH} that some operators on the Sweedler  algebra $H_4$, satisfying only the condition (\ref{0}),  are RB-operators on the adjoint Lie algebra $H^{(-)}_4$.
It is natural to formulate   the next definitions. Let $A$ be an associative algebra.  We shall call Rota--Baxter operators on $A$  \emph{associative RB-operators} (A-RB-operators, shortly). On the vector space $A$  we can construct the  adjoint Lie algebra $A^{(-)}$ and find its Rota--Baxter operators on this algebra. We shall call these operators {\it Lie RB-operators} (L-RB-operators, shortly). Also, we can construct the Jordan algebra  $A^{(+)}$ and find its Rota--Baxter operators, which we shall call {\it Jordan  RB-operators} (J-RB-operators, shortly).  It is interesting to find connections between these operators. It is easy to prove that any associative  RB-operator on $A$ is a Lie  RB-operator on $A^{(-)}$ and is a Jordan  RB-operator on $A^{(+)}$.

 In this paper we describe all RB-operators on the adjoint Lie algebra $H^{(-)}_4$.

The paper is organised as follows.
In~\S~\ref{Prel}, we state the required preliminaries on Hopf algebras, Rota--Baxter operators on algebras. In Section \ref{Ker-2}, we found the RB-operators on the three-dimensional Lie algebras, which are the subalgebras of $H^{(-)}_4$.
Further we  describe RB-operators on $H^{(-)}_4$ what have the kernel of  dimension  3. In Sections 4 and 5 we describe Roth--Baxter operators with the kernels of dimensions 2 and 1, respectively.
Finally, the non-degenerate RB-operators on $H^{(-)}_4$ are described in  Section 5.

\bigskip


\section{Preliminaries} \label{Prel}

\subsection{Hopf algebras} Recall some facts on  Hopf algebras (see, for example, \cite[Chapter~11]{T}).

 Let $A$ be an associative algebra with unit $1_A$ over a commutative ring with unit, $\Bbbk$. Assume that $A$ is provided with multiplicative
$\Bbbk$-linear homomorphisms
$$
\Delta \colon A \to A^{\otimes 2} = A \otimes_{\Bbbk} A~\mbox{and}~\varepsilon \colon A \to \Bbbk,
$$
called
the comultiplication and the counit respectively, and a $\Bbbk$-linear map $S \colon A \to A$, called the antipode. It is understood that $\Delta(1_A) = 1_A \otimes 1_A$ and $\varepsilon(1_A) = 1$. The tuple $(A, \Delta, \varepsilon, S)$ is said to be a Hopf  algebra if these
homomorphisms satisfy together with the multiplication $m \colon A \times A \to A$ the following identities:
\begin{equation}\label{1.1.a}
(\id_A  \otimes \Delta) \Delta  = (\Delta \otimes \id_A) \Delta,
\end{equation}
\begin{equation}\label{1.1.b}
m (S \otimes \id_A ) \Delta  = m (\id_A \otimes S) \Delta = \varepsilon \cdot 1_A,
\end{equation}
\begin{equation}\label{1.1.c}
(\varepsilon \otimes \id_A ) \Delta  = (\id_A \otimes  \varepsilon) \Delta = \id_A.
\end{equation}
Note that to  write down the first equality we identify $(A \otimes A) \otimes A = A \otimes (A \otimes A)$ via $(a \otimes b) \otimes c = a \otimes (b \otimes c)$, where $a, b, c \in A$.  Similarly, to write down the third equality, we identify $\Bbbk \otimes A = A \otimes \Bbbk = A$ via $1 \otimes a = a \otimes 1 = a$. The axioms imply that the antipode $S$ is an antiautomorphism of both the
algebra and the coalgebra structures in~$A$. This means that
$$
m (S \otimes S) = S \circ m \circ P_A \colon A^{\otimes 2} \to A,~~~~P_A (S \otimes S) \Delta = \Delta  \circ S  \colon A \to A^{\otimes 2},
$$
where $P_A$ denotes the flip
$$
a \otimes b \mapsto b \otimes a \colon A^{\otimes 2} \to A^{\otimes 2}.
$$
 It also follows from the axioms that $S(1_A) = 1_A$ and $\varepsilon \circ S = \varepsilon \colon A \to \Bbbk$.

\begin{example}
 The famous non-cocommutative 4-dimension Sweedler  algebra $H_4$ as algebra is generated by two elements $x$, $g$ with multiplication
$$
g^2 = 1,~~x^2 = 0,~~x g = -g x.
$$
Comultiplication, counit and antipod are defined by the rules
$$
\Delta(g) = g \otimes g,~~\Delta(x) = x \otimes 1 + g \otimes x,~~ \Delta(g x) = g x \otimes g + 1 \otimes g x,
$$
$$
\varepsilon(g) = 1,~~ \varepsilon(x) = 0,~~ \varepsilon(gx) = 0,
$$
$$
 S(g) = g^{-1} = g,~~S(x) = -g x ,~~ S(g x) = x .
$$
The antipode $S$ has order 4 and for any $a \in H_4$ we have $S^2(a) = g a g^{-1}$.
\end{example}

 Throughout the paper  we shall use the
Sweedler's notation with superscripts for all the comultiplications,
$$
\Delta(h) = h_{(1)} \otimes h_{(2)},~~~h \in H.
$$

\subsection{Rota--Baxter operators}
As has already been noted for its popularity Rota--Baxter operators are obliged of the work G.~Baxter~\cite{Baxter}.
For basic results and the main properties of Rota--Baxter algebras see~\cite{GuoMonograph}.

\begin{definition}
Let $A$ be an algebra over a field~$\Bbbk$. A linear operator $R$ on $A$ is called a Rota--Baxter operator of weight~$\lambda\in \Bbbk$ if
\begin{equation}\label{RBAlgebra}
R(a)R(b) = R\left( R(a) b + a R(b) + \lambda a b \right)
\end{equation}
for all $a, b \in A$.
An algebra endowed with a Rota--Baxter operator is called a~Rota--Baxter algebra.
\end{definition}

In \cite{Ma} (see also \cite{Ma-1}) were found algebra RB-operators on $H_4$ as on associative algebra. It was proven that any such operator has one of the following form:\\

\begin{itemize}
\item[(a)] $R(1) = 0$,  $R(g) = 0$, $R(x) = -\lambda x$, $R(gx) = -\lambda gx$;

\item[(b)] $R(1) = -\lambda 1$,  $R(g) = -\lambda g$, $R(x) = 0$, $R(gx) = 0$;

\item[(c)] $R(1) = -\lambda 1$,  $R(g) = -\lambda g$, $R(x) = -\lambda x$, $R(gx) = -\lambda gx$;

\item[(d)] $R(1) = 0$, $R(g) = -p_1 1 + p_1 g - \frac{(\lambda + p_1) (\lambda + p_1 + p_2)}{p_3} x + \frac{(\lambda + p_1) (\lambda + p_2)}{p_3} gx$,\\
 $R(x) = -p_3 1 + p_3 g - (2 \lambda + p_1 + p_3) x +  (\lambda + p_2) g x$,\\
   $R(g x) = -p_3 1 + p_3 g - (\lambda + p_1 + p_2) x +  p_2 g x$;

\item[(e)] $R(1) = -\lambda 1$, $R(g) = (\lambda + p_1) 1 + p_1 g - \frac{(\lambda + p_1) (\lambda + p_1 + p_2)}{p_3} x + \frac{(\lambda + p_1) (\lambda + p_2)}{p_3} gx$,\\
 $R(x) = p_3 1 + p_3 g - (2 \lambda + p_1 + p_2) x +  (\lambda + p_2) g x$,\\
   $R(g x) = p_3 1 + p_3 g - (\lambda + p_1 + p_2) x +  p_2 g x$;

\item[(f)] $R(1) = -\lambda 1$, $R(g) = \lambda  1 + p_1  x + \frac{ p_1 p_2}{\lambda + p_2} gx$,\\
$R(x) = - (\lambda +  p_2) x - p_2 g x$, $R(gx) = (\lambda +  p_2) x + p_2 g x$;

\item[(g)] $R(1) = -\lambda 1$, $R(g) = \lambda 1  + \frac{\lambda (\lambda + p_1)}{p_2} x + \frac{\lambda (\lambda + p_1)}{p_2} gx$,\\
 $R(x) = -p_2 1 - p_2 g - (2 \lambda + p_1) x -  (\lambda + p_1) g x$,\\
   $R(g x) = p_2 1 + p_2 g + (\lambda + p_1) x +  p_1 g x$;

\item[(h)]   $R(1) = \frac{1}{2} \lambda  1 -  \frac{1}{2} \lambda g  + p_1 x +  p_2 g x$, $R(g) = \frac{1}{2} \lambda  1 -  \frac{1}{2} \lambda g - p_2 x +  p_1 g x$,\\
$R(x) = -\frac{1}{2} \lambda  x -  \frac{1}{2} \lambda g x$, $R(g x) = -\frac{1}{2} \lambda  x -  \frac{1}{2} \lambda g x$.
\end{itemize}

\bigskip


\section{  Rota--Baxter operators on subalgebras of  Lie algebra   $H_4^{(-)}$} \label{Ker-2}

Let us prove a general fact on  RB-operators of different types.

\begin{proposition}
Let $R$ be an A-RB-operator of weight $\lambda$ on an associative algebra $A$ and $A^{(-)}$ is its adjoint  Lie algebra. Then
$R$ is an L-RB-operator of weight $\lambda$ on the Lie algebra   $A^{(-)}$ and $R$ is an J-RB-operator of weight $\lambda$ on the Jordan algebra   $A^{(+)}$.
\end{proposition}

\begin{proof}
For any $a, b \in A$ hold
$$
R(a) R(b) = R(R(a) b + a R(b) + \lambda a b),~~R(b) R(a) = R(R(b) a + b R(a) + \lambda b a).
$$
 Subtract the second equality from the first,
$$
[R(a),  R(b)] = R ( [R(a),  b] + [a,  R(b)] + \lambda [a,  b]).
$$
Hence, $R$ is an  RB-operator of weight $\lambda$ on the Lie algebra   $A^{(-)}$.

A proof for $A^{(+)}$ is similar.

\end{proof}

The following properties of $H_4^{(-)}$
are well known.

1) $[H_4^{(-)}, H_4^{(-)}] = \Bbbk \cdot  x + \Bbbk \cdot  gx$;

2) $[[H_4^{(-)}, H_4^{(-)}], [H_4^{(-)}, H_4^{(-)}]]$ = 0, i. e. $H_4^{(-)}$ is an metabelian algebra;

3) subalgebra $\Bbbk \cdot 1$ lies in the center of $H_4^{(-)}$;

4) the set $I = \Bbbk \cdot 1 + \Bbbk \cdot  x + \Bbbk \cdot  gx$ is an ideal in $H_4^{(-)}$ and $[I, I] = 0$.

\medskip

In the sequel, we will assume that  $char(\Bbbk) \neq 2$ and  will denote  by letters of the Greek alphabet, except $\phi$ and $\psi$,  elements of the field $\Bbbk$.

If we put  $e=x+gx, f=x-gx$, then
$$
[1,g]=[1,e]=[1,f] =0,~~[g,e]=2e,~~ [g,f]=-2f,~~ [e,f]=0.
$$

The map $\phi  \colon H_4^{(-)} \to H_4^{(-)},$ which is defined by the rule
$$
\phi(1)=1,~~\phi(g)=-g,~~\phi(e)=-f, ~~\phi(f)=-e
$$ is an automorphism of the Lie algebra $H_4^{(-)}$ and the Hopf algebra  $H_4$.

The vector spaces  $K=\Bbbk e+\Bbbk f$,  $K_e=\Bbbk 1+ \Bbbk e$ and  $K_f= \Bbbk 1+ \Bbbk f$ are two dimensional ideals of the algebra $H_4^{(-)}$.

\begin{lemma} \label{lm0}
    Let  $L$ be a two dimensional ideal of the  algebra  $H_4^{(-)}$, which is different from   $K$. Then $L=K_e$ or $L=K_f$.
\end{lemma}

\begin{proof}
Let $L\neq K$. Then  $L$ contains a a non-zero element
$u=\alpha 1+\beta e+\gamma f$ such that $\alpha\neq 0$ and we have
$$
[g,u]=2\beta e-2\gamma f\in L, ~~ [g, [g,u]]=2\beta e+2\gamma f  \in L.
$$
Therefore,  $\alpha 1,\beta e,\gamma f\in L$.
Consequently,  $L=K_e$ or $L=K_f$.
\end{proof}

The vector spaces $I=\Bbbk 1+ \Bbbk e+\Bbbk  f$ and  $J_\alpha=\Bbbk(\alpha 1+g)+\Bbbk e+\Bbbk f$, where
$\alpha\in \Bbbk$  are  three-dimensional ideals of   $H_4^{(-)}$. Put $J=J_0$.

\begin{lemma} \label{lm1} Let $L$ be a three-dimensional non-abelean subalgebra of $H_4^{(-)}$.
Then $L=span(\alpha 1+g,e,f)$ or $L=span(1, g+\gamma e,f)$ or
 $L=span(1, g+\gamma f,e) $. If further $L$ is a three-dimensional ideal of  $H_4^{(-)}$, then  $L=I$ or $L=J_\alpha$.
\end{lemma}

\begin{proof}  Let $L$ be a three  dimensional non-abelean subalgebra.
Since   $\dim (L\cap I)\geq 2$, then there are  $a,b\in L\cap I$, which
are linear  independent. Let  $c=\alpha 1+\beta g+\gamma e+\delta f$
be an  element of   $L$ such that $a$, $b$, $c$ is a basis of $L$. Since
$L$ is not abelean,  we can assume that  $\beta=1$ and $[c,a]\neq 0$ or $[c, b]\neq 0$.

 Let  $a=\xi 1+\mu e+\eta f$ and  $[a, c]\neq 0$. Then
 $2\mu e-2\eta f=[c,a]$ and  $\mu\neq 0$ or $\eta \neq 0$.
Since   $\mu e-\eta f\in L$, we have $$2\mu e+2\eta f=[c,\mu e-\eta f]\in L.$$ Therefore, $ e\in L$ or $f\in L$.

If  $e, f\in  L$ then  $L=span(\alpha 1+ g,e,f)$.

Assume  that $e\not\in L$. Then $f\in L$ and we may suppose  that  $a=\xi 1+\mu e.$ Since
$a,b$ are linear  independent, we can  suppose that $b= 1.$
Hence, $L=span(1,g+\gamma e,f)$.

Similarly $L=span(1,g+\gamma f,e)$, if $e\in L$.

Let $L$ be an  ideal of  $H_4^{(-)}$. If  $L$  is not abelean then according
to the fact which was proven above we have  $L=J_\alpha$ for some  $\alpha\in \Bbbk$.

 Let   $L$ be an abelean ideal and let  $u=\alpha1+\beta g+\gamma  e+\delta f\in L$.
In this case   $2\beta e=[u,e]\in L$. Since $L$ is abelean,  $2\beta e=[u, e]  =0$. Consequently,   $L=I$.
\end{proof}

Recall a definition of a Rota-Baxter operator on a Lie algebra. A linear map $R \colon H_4^{(-)} \to H_4^{(-)}$ is called a Rota-Baxter
operator (RB-operator) of weight $\lambda \in \Bbbk$, if the following
condition holds
\begin{equation}\label{YBop}[R(v),R(u)]=R([R(v),u]+[v,R(u)]+\lambda[u,v]),~~~v, u \in H_4^{(-)}.
\end{equation}

It is clearly that the image  ${\rm Im}R$ of $R$ is a subalgebra of  $H_4^{(-)}$.

For any
$a\in H_4^{(-)}$ and $v\in \ker R$ we have
$$
0= [R(a),R(v)]=R([R(a),v]+\lambda [a, v]).
$$
 Therefore,
\begin{equation}\label{kercon}
[R(a)+\lambda a,v]\in \ker R, \,\,  a\in H_4^{(-)},~~v  \in \ker R.
\end{equation}
Moreover, we see that $[R(a)+\lambda a,v]\in  \Bbbk e+ \Bbbk f.$
Thus, \begin{equation}\label{kercon1}[R(a)+\lambda a,v]\in \ker R\cap K.\end{equation}

\medskip

Our  aim  to describe Rota--Baxter operator of non-zero weight on the algebra   $H_4^{(-)}$.

Let  $\lambda\neq 0$. In this case the kernel $\ker R$ of operator $R$ is a subalgebra of $H_4^{(-)}$.
Indeed, for any   $v,u\in \ker R$ we have
$$
\lambda R([v,u])=[R(v),R(u)]-R([R(v),u]-[v,R(u)])=0.
$$

The following   lemmas  describe Rota--Baxter operators on some three-dimensional
  non-abelean Lie algebra  which is  a subalgebra of $H_4^{(-)}$.

\begin{lemma} \label{lm2} Let  $L$  be a three-dimensional  Lie algebra with a basis $h,e,f$ and  the following
 multiplication table
    $$[h,e]=2e,\,[h,f]=-2f,\, [e,f]=0.$$  Then any  Rota--Baxter  operator $R$ on  $L$  has one of the following form

$ 1.1.\, R(h)=\alpha_2^{-1}(\beta_2+\lambda)R(e),~~  R(e)=\alpha_2h+\beta_2e+\gamma_2f, ~~R(f)=0.$

$ 1.2.\, R(h)=-\alpha_2^{-1}((\beta_2+\lambda) R(e)+\gamma_2R(f)),~~ R(e)=\alpha_2h+\beta_2e+\gamma_2f, ~~R(f)=-\lambda f.$

$1.3.\, R(h)=\alpha_1h+\beta_1e+\gamma_1f,~~R(e)=R(f)=0.$

$1.4.\, R(h)=-\frac{\lambda}{2}h+\beta_1e+\gamma_1f, ~~R(e)=0, ~~R(f)=\beta_3e,\beta_3\neq 0.$

$1.5.\, R(h)=\alpha_1h+\beta_1e+\gamma_1f, ~~R(e)=0, ~~R(f)=-\lambda f.$

$1.6.\, R(h)=\beta_1e+\gamma_1f,~~R(e)=0,R(f)=\beta_3e-\lambda f,~~\beta_3\neq 0. $


$1.7.\, R(h)=-\lambda h+\beta_1e+\gamma_1f,~~R(e)=-\lambda e, ~~R(f)=\beta_3 e,\beta_3\neq 0.$

$1.8.\, R(h)=\alpha_1h+\beta_1e+\gamma_1f,~~ R(e)=-\lambda e, ~~R(f)=-\lambda f.$

$1.9.\, R(h)=-\frac{\lambda}{2}h+\beta_1e+\gamma_1f,~~R(e)=-\lambda e,~~ R(f)=\beta_3e-\lambda f,\beta_3\neq 0.$


$1.10. \, R(h)=-\frac{\lambda }{2}h+\beta_1e+\gamma_1f,~~R(e)=\beta_2e +\gamma_2f,~~R(f)=\gamma_2^{-1}\beta_2(\beta_2+\lambda)e+\beta_2 f, $

$\beta_2\neq 0, -\lambda, \gamma_2 \neq 0. $
\end{lemma}
\begin{proof} Note that the map $\psi \colon L \to L$, which is defined by the rule
	$$\psi(h)=-h,\psi(e)=-f,\psi(f)=-e$$ is an automorphism of  $L$.
	
	Let
	$$R(h)=\alpha_1h+\beta_1e+\gamma_1
	f,R(e)=\alpha_2h+\beta_2e+\gamma_2f,$$
	$$R(f)=\alpha_3h+\beta_3e+\gamma_3f.$$
	Applying (\ref{YBop}) for $e, f$ we obtain
	\begin{equation}\label{e-f}\left\{ \begin{array}{l}
			\alpha_2\alpha_3=0, \\
			\alpha_2\beta_3=0, \\
			\alpha_3\gamma_2=0.
		\end{array}\right.
	\end{equation}
	
	Applying (\ref{YBop}) for $h,e $ we get
	\begin{equation}\label{h-e}\left\{ \begin{array}{l}
			\alpha_1\beta_2-\alpha_2\beta_1=(\beta_2+\lambda+\alpha_1)\beta_2-\gamma_2\beta_3, \\
			\alpha_2\gamma_1-\alpha_1\gamma_2=(\beta_2+\lambda+\alpha_1)\gamma_2-\gamma_2\gamma_3, \\
			(\beta_2+\lambda+\alpha_1)\alpha_2=\gamma_2\alpha_3.
		\end{array}\right.
	\end{equation}
	Similar for  $h,f$,
	\begin{equation}\label{h-f}
		\left\{ \begin{array}{l}
			\alpha_3\beta_1-\alpha_1\beta_3=(\gamma_3+\lambda+\alpha_1)\beta_3-\beta_2\beta_3, \\
			\alpha_1\gamma_3-\alpha_3\gamma_1=(\gamma_3+\lambda+\alpha_1)\gamma_3-
			\beta_3\gamma_2, \\
			(\gamma_3+\lambda+\alpha_1)\alpha_3=\alpha_2\beta_3.
		\end{array}\right. \end{equation}

	Let $\alpha_2\neq 0$. Then from (\ref{e-f}) $\alpha_3=\beta_3=0$. From (\ref{h-e})
	$\alpha_1=-\beta_2-\lambda$, $\beta_1=-\alpha_2^{-1}(\beta_2+\lambda)\beta_2$,
	$\gamma_1=-\alpha_2^{-1}(\beta_2+\lambda+\gamma_3)\gamma_2$.  From  (\ref{h-f}) $\gamma_3(\gamma_3+\lambda)=0$.
	Thus, we have the case 1.1:
	$$ R(h)=-\alpha_2^{-1}(\beta_2+\lambda)R(e),~~R(e)=\alpha_2h+\beta_2e+\gamma_2f, ~~R(f)=0,
	$$ or the case   1.2:
	$$
	R(h)=-\alpha_2^{-1}((\beta_2+\lambda)R(e)+\gamma_2R(f)),~~
	R(e)=\alpha_2h+\beta_2e+\gamma_2f, R(f)=-\lambda f.$$
	
	Suppose that $\alpha_2=0$. If $\alpha_3\neq0$, then from (\ref{e-f}) $\gamma_2=0$, from (\ref{h-e}) $\beta_2(\beta_2+\lambda)=0$, from (\ref{h-f}) $\alpha_1=-\gamma_3-\lambda$, $\beta_1=-\alpha_3^{-1}(\gamma_3+\lambda+\beta_2)\beta_3$ and $\gamma_1=-\alpha_3^{-1}(\gamma_3+\lambda)\gamma_3$.
	Thus,
	$$ R(h)=-\alpha_3^{-1}(\gamma_3+\lambda)R(f), ~~R(e)=0, ~~R(f)=\alpha_3h+\beta_3e+\gamma_3f
	$$
	or
	$$
	R(h)=-\alpha_3^{-1}((\gamma_3+\lambda)R(f)+\beta_3R(e)),~~R(e)=-\lambda e,~~ R(f)=\alpha_3h+\beta_3e+\gamma_3f.
	$$
	Up   to automorphism    $\psi$ these  two cases are coincided     with  the cases 1.1 and  1.2.
	
	Suppose that $\alpha_2=\alpha_3=0$. Then from  (\ref{h-e}) and  (\ref{h-f})
	we obtain
	\begin{equation}\label{h-e-f}
		\left\{ \begin{array}{l}
			(\beta_2+\lambda)\beta_2=\gamma_2\beta_3,\\
			(\beta_2+\lambda+2\alpha_1-\gamma_3)\gamma_2=0,\\
			(\gamma_3+\lambda)\gamma_3=\gamma_2\beta_3.\\
			(\gamma_3+\lambda+2\alpha_1-\beta_2)\beta_3=0.
		\end{array}\right. \end{equation}

	Let $\beta_2=\gamma_2=0$. Then from  (\ref{h-e-f}) \begin{equation}\nonumber
		\left\{ \begin{array}{l}
			(\gamma_3+\lambda)\gamma_3=0,\\
			(\gamma_3+\lambda+2\alpha_1)\beta_3=0.
		\end{array}\right. \end{equation}
	Thus, we have the case  1.3:
	$$ R(h)=\alpha_1h+\beta_1e+\gamma_1f,~~R(e)=R(f)=0,$$
	or the case  1.4:
	$$R(h)=-\frac{\lambda}{2}h+\beta_1e+\gamma_1f, ~~R(e)=0, ~~R(f)=\beta_3e,\beta_3\neq 0,$$
	or the case  1.5:
	$$
	R(h)=\alpha_1h+\beta_1e+\gamma_1f,~~R(e)=0, ~~R(f)=-\lambda f,
	$$
	or the case  1.6:
	$$
	R(h)=\beta_1e+\gamma_1f,~~R(e)=0,~~R(f)=\beta_3e-\lambda f,\beta_3\neq 0.
	$$
	
	Let $\beta_2\neq 0,\gamma_2=0$. Then from (\ref{h-e-f}) \begin{equation}\nonumber
		\left\{ \begin{array}{l}
			\beta_2=-\lambda,\\
			(\gamma_3+\lambda)\gamma_3=0,\\
			(\gamma_3+2\lambda+2\alpha_1)\beta_3=0.
		\end{array}\right. \end{equation}
	Thus, we have the case
	$$
	R(h)=\alpha_1h+\beta_1e+\gamma_1f,~~R(e)=-\lambda e,~~R(f)=0,
	$$
 which coincides up to automorphism $\psi$ with the case 1.5,
 or the case  1.7:
	$$
	R(h)=-\lambda h+\beta_1e+\gamma_1f,~~R(e)=-\lambda e, ~~R(f)=\beta_3 e,~~\beta_3\neq, 0$$ or the case 1.8:
	$$R(h)=\alpha_1h+\beta_1e+\gamma_1f, R(e)=-\lambda e, R(f)=-\lambda f$$ or the case  1.9:
	$$R(h)=-\frac{\lambda}{2}h+\beta_1e+\gamma_1f,R(e)=-\lambda e, R(f)=\beta_3e-\lambda f,\beta_3\neq 0.$$
	
		Let $\beta_2=0$, $\gamma_2\neq 0$. Then from  (\ref{h-e-f})
	\begin{equation}\nonumber
		\left\{ \begin{array}{l}
			\beta_3=0,\\
			\gamma_3=\lambda+2\alpha_1,\\
			(\gamma_3+\lambda)\gamma_3=0.
		\end{array}\right. \end{equation}
	
	Thus, we have the case
	$$
	R(h)=-\frac{\lambda }{2}h+\beta_1e+\gamma_1f,~~R(e)=\gamma_2f,~~R(f)=0, \gamma_2\neq 0,
      $$
 which coincides up to automorphism $\psi$ with the case 1.4
	or the case
	$$
	R(h)=-\lambda h+\beta_1e+\gamma_1f,~~R(e)=\gamma_2f,~~R(f)= -\lambda f,~~\gamma_2\neq 0,
	$$ which coincides   up   to autoomorphism   $\psi$ with the case 1.7.
	
	Let $\beta_2\neq 0$, $\gamma_2\neq 0$. Then from (\ref{h-e-f})
	\begin{equation}\nonumber
		\left\{ \begin{array}{l}
			(\beta_2+\lambda)\beta_2=\gamma_2\beta_3,\\
			\gamma_3=\beta_2+\lambda+2\alpha_1,\\
			(\gamma_3+\lambda)\gamma_3=\gamma_2\beta_3,\\
			(\lambda+2\alpha_1)\beta_3=0.
		\end{array}\right. \end{equation}
	Thus, we have the case
	$$
	R(h)=\beta_1e+\gamma_1f,~~R(e)=-\lambda e+\gamma_2 f,~~R(f)=0,~~\gamma_2\neq 0,
	$$
	which coincides up to automorphism   $\psi$ with the case 1.6  or  the case
	$$ R(h)=-\frac{\lambda}{2} h+\beta_1e+\gamma_1f,~~R(e)=-\lambda e+\gamma_2f,~~R(f)= -\lambda f,~~\gamma_2\neq 0,
	$$
	which coincides up to automorphism  $\psi$ with the case 1.9,   or  the case  1.10:
	$$
	R(h)=-\frac{\lambda }{2}h+\beta_1e+\gamma_1f,~~R(e)=\beta_2e +\gamma_2f,
	$$
	$$R(f)=\gamma_2^{-1}\beta_2(\beta_2+\lambda)e+\beta_2 f,  ~~\beta_2\neq 0, ~~\lambda, \gamma_2 \neq 0.
	$$
\end{proof}

\begin{lemma} \label{lm3} Let  $L$ be a three-dimensional  Lie algebra with a basis  $h$, $y$, $z$ and
	the following multiplication table
	$$
	[h,y]=2 y,~~[h,z]=[y,z]=0.
	$$
	Then each Rota-Baxter operator  $R$ on   $L$ has one of the  following forms:
	
	$2.1.\,\, R(h)=\alpha_1 h+\beta_1 y+\gamma_1 z,~~ \alpha_1\neq 0, -\lambda,~~R(y)=0,~~R(z)=\alpha_3h+\alpha_1^{-1}\beta_1\alpha_3y+\gamma_3z.$
	
	$2.2.\,\, R(h)=\beta_1 y+\gamma_1z, R(y)=0, R(z)=\beta_3y+\gamma_3 z.$
	
	$2.3.\,\, R(h)=\gamma_1z,~~R(y)=0, R(z)=\alpha_3h+\beta_3 y+\gamma_3 z.$
	
	$2.4.\,\, R(h)=\alpha_1h+\beta_1 y+\gamma_1z,~~R(y)=-\lambda y,~~ R(z)=\gamma_3 z,\alpha_1\neq 0.$
	
	$2.5.\,\, R(h)=-\lambda h+\beta_1 y+\gamma_1z,~~ R(y)=-\lambda y, ~~R(z)=\beta_3 y+\gamma_3 z.$
	
	$2.6.\,\, R(h)=\alpha_1h+\beta_3\alpha_3^{-1}(\alpha_1+\lambda) y+\gamma_1z,~~R(y)=-\lambda y,~~
	R(z)=\alpha_3h+\beta_3 y+\gamma_3 z,\alpha_1\neq 0.$
	
	$2.7.\,\, R(h)=\beta_1 y+\gamma_1z, R(y)=-\lambda y+\gamma_2z,~~R(z)=\alpha_3h+\lambda^{-1}\alpha_3\beta_1y+\gamma_3z.$
	
	$2.8.\,\, R(h)=-\lambda h+\beta_1 y+\gamma_1z,~~ R(y)=\gamma_2z,~~R(z)=\alpha_3h-\lambda^{-1}\alpha_3\beta_1y+\gamma_3z.$
	
	$ 2.9.\,\, R(h)=\alpha_1h-\alpha_2^{-1}\alpha_1(\alpha_1+\lambda)y+\gamma_1z,~~R(y)=\alpha_2h-(\alpha_1+\lambda)y+\gamma_2 z, R(z)=\gamma_3z.$
\end{lemma}

\begin{proof}
	Let	$$R(h)=\alpha_1h+\beta_1y+\gamma_1z,~~R(y)=\alpha_2h+\beta_2y+\gamma_2z,~~R(z)=\alpha_3h+\beta_3y+\gamma_3z.$$
	Applying  (\ref{YBop}) for $h$, $y$ we obtain
	
	\begin{equation}\label{x-y}\left\{ \begin{array}{l}
			(\alpha_1+\beta_2+\lambda)\alpha_2=0,\\
			\alpha_1\beta_2-\alpha_2\beta_1=(\alpha_1+\beta_2+\lambda)\beta_2,\\
			(\alpha_1+\beta_2+\lambda)\gamma_2=0.
		\end{array}\right.
	\end{equation}
	Applying  (\ref{YBop}) for $h$, $z $ we obtain
	\begin{equation}\label{x-z}\left\{ \begin{array}{l}
			\beta_3\alpha_2=0,\\
			\alpha_1\beta_3-\alpha_3\beta_1=\beta_3\beta_2,\\
			\beta_3\gamma_2=0.
		\end{array}\right.
	\end{equation}
	Similar for  $y$, $z$
	\begin{equation}\label{y-z}
		\left\{ \begin{array}{l}
			\alpha_2\beta_3=0,\\
			\alpha_3\alpha_2=0,\\
			\alpha_3\gamma_2=0.
		\end{array}\right. \end{equation}
	
	Suppose $\alpha_1+\beta_2+\lambda\neq 0.$ Then from (\ref{x-y})
	$\alpha_2=\gamma_2=0$ and $(\beta_2+\lambda)\beta_2=0$.
	
	We assume  $\beta_2=0$. Then $\alpha_1\neq -\lambda$ and from   (\ref{x-z}) $\alpha_1\beta_3=\alpha_3\beta_1$. If $\alpha_1\neq 0$, then
	$$R(h)=\alpha_1h+\beta_1y+\gamma_1 z,~~R(y)=0,~~R(z)=\alpha_3h+\alpha_1^{-1}\beta_1\alpha_3y+\gamma_3z.
	$$
	If  $\alpha_1= 0$,  then $\alpha_3\beta_1=0$.  Thus we obtain
	$$R(h)=\beta_1 y+\gamma_1z, ~~R(y)=0, ~~R(z)=\beta_3y+\gamma_3 z,$$
	or
	$$
	R(h)=\gamma_1z,~~ R(y)=0, ~~R(z)=\alpha_3h+\beta_3 y+\gamma_3 z.
	$$
	If $ \beta_2=-\lambda$, then  $\alpha_1\neq 0$ and from (\ref{x-z})
	$\beta_3(\alpha_1+\lambda)=\alpha_3\beta_1$. Let
	$\alpha_3=0$. Then  $\beta_3(\alpha_1+\lambda)=0$.
	Thus we have
	$$
	R(h)=\alpha_1h+\beta_1 y+\gamma_1z,~~R(y)=-\lambda y,~~ R(z)=\gamma_3 z,\alpha_1\neq 0,
	$$
	or
	$$
	R(h)=-\lambda h+\beta_1 y+\gamma_1z,~~R(y)=-\lambda y,~~ R(z)=\beta_3 y+\gamma_3 z.
	$$

	If  $\alpha_3\neq 0$, then
	$$R(h)=\alpha_1h+\beta_3\alpha_3^{-1}(\alpha_1+\lambda) y+\gamma_1z,~~ R(y)=-\lambda y,~~ R(z)=\alpha_3h+\beta_3 y+\gamma_3z,~~\alpha_1\neq 0.
	$$

	Suppose $\alpha_1+\beta_2+\lambda=0.$  From (\ref{x-y}) follows
	$\alpha_1\beta_2=\alpha_2\beta_1$.
	
	We assume that  $\alpha_2=0$. Then  $\alpha_1\beta_2=0$. If   $\alpha_1=0$, then  $\beta_2=-\lambda$ and  from  (\ref{x-z}) $\beta_3=
	\lambda^{-1}\alpha_3\beta_1$. Thus
	$$
	R(h)=\beta_1 y+\gamma_1z, ~~R(y)=-\lambda y+\gamma_2z,~~R(z)=\alpha_3h+\lambda^{-1}\alpha_3\beta_1y+\gamma_3z.
	$$
	If $\alpha_1\neq 0$, then  $\beta_2=0, \alpha_1=-\lambda$ and  from  (\ref{x-z}) follows $\beta_3=- \lambda^{-1}\alpha_3\beta_1$. Thus
	$$
	R(h)=-\lambda h+\beta_1 y+\gamma_1z, R(y)=\gamma_2z,~~R(z)=\alpha_3h-\lambda^{-1}\alpha_3\beta_1y+\gamma_3z.
	$$
	
	Let $\alpha_2\neq 0$. Then from (\ref{x-z}) and  (\ref{y-z}) follows $\alpha_3=\beta_3=0$. Thus
	$$
	R(h)=\alpha_1h-\alpha_2^{-1}\alpha_1(\alpha_1+\lambda)y+\gamma_1 z,
	~~R(y)=\alpha_2h-(\alpha_1+\lambda)y+\gamma_2 z, ~~R(z)=\gamma_3z.
	$$
\end{proof}

Let us consider the case  when the kernel  $\ker R$ has the dimension 3. By Lemma~\ref{lm2} and Lemma~\ref{lm3} we may suppose that  $\ker R$ has
the basis $\alpha1+g ,e,f$ or the basis $1,g+\gamma f,e$.

Let  $\alpha1+g ,e,f$  be a basis of $\ker R$ and let
$$
R(1)=\alpha_1 1+\alpha_2 g+\alpha_3 e+\alpha_4 f.
$$
Then  $\alpha_i \neq 0$ for some  $\alpha_i$, $R(g)=-\alpha R(1)$, $R(e)=R(f)=0$.

Let  $1$, $g+\gamma f$, $e$ be a  basis of   $\ker R$ and let
$$
R(f)=\alpha_1 1+\alpha_2 g+\alpha_3 e+\alpha_4 f.
$$
Then $R(f)\neq 0$ and $R(g)=- \gamma R(f)$.  By (\ref{YBop}), we have
$$
0=[R(g),R(f)]=R([R(g),f]+[g,R(f)]-2\lambda f)=2(\gamma\alpha_2 -\alpha_4-\lambda)R(f).
$$ Hence $\alpha_4=\gamma\alpha_2-\lambda$.

We have proven the next result.

\begin{theorem}
	Let $\dim\ker R=3$.  Then each Rota--Baxter operator $R$ on the algebra $H_4^{(-)}$ has a form
	$$
	R(1)=\alpha_1 1+\alpha_2 g+\alpha_3 e+\alpha_4 f, ~~R(g)=-\alpha R(1), ~~R(e)=R(f)=0,
	$$
	where one  from elements  $\alpha_1,\alpha_2,\alpha_3,\alpha_4$ is non-zero, or a form
	$$
	R(1)=0, ~~R(g)=-\gamma R(f), ~~R(e)=0, ~~R(f)=\alpha_1 1+\alpha_2 g+\alpha_3 e+(\gamma\alpha_2-\lambda)f,
	$$
	where  one from elements  $\alpha_1,\alpha_2,\alpha_3$ is non-zero or
	$\gamma\alpha_2\neq \lambda$.
\end{theorem}

\bigskip


\section{The kernel $R$ has  the dimension two}

Suppose that  $\ker R$ has the  dimension 2. In this case ${\rm Im}  R$ also has the dimension 2 and  we have that
$\ker R\cap I\neq 0$ and ${\rm Im} R\cap I\neq 0$.

\subsection{The kernel $\ker R$ is abelean}

Let  $[\ker R,\ker R]=0$. Assume
$u=\alpha_0 1+\beta_0 e+\gamma_0
f$ is a non-zero element of  $\ker R\cap I$ and  $v=\delta_0 g+w\in \ker R$, $w\in
I$, $\delta_0\in \Bbbk$.

\begin{lemma}  Let  $\delta_0\neq 0$. Then  $v,u$  is a basis of  $\ker R$ and  $R$ has one of the following forms:
	
	{\rm (i)}   $$R(1)=0, ~~R(g)=0, ~~R(e)=\alpha 1- \lambda e,~~R(f)=\alpha_1 1-\lambda f.$$
	
	{\rm (ii)} $$R(1)=0,~~ R(g)=0, ~~R(e)=\alpha 1+ \sigma g-\lambda e,~~R(f)=-\lambda f.$$
	
	{\rm (iii)}
	$$R(1)=0, ~~R(e)=-\lambda e+\alpha 1,~~R(f)=-\lambda f+\alpha_1 1,$$
	$$ R(g)=\xi\lambda e+\eta \lambda f-(\xi\alpha+\eta\alpha_1)1, ~~\xi\neq 0 ~\text{ or }~ \nu\neq 0.$$
	
	{\rm (iv)}  $$
	R(e)=-\lambda e,~~R(f)=-\lambda f+\alpha_2 1+\beta_2 (g+\xi e+\eta f).$$
	$$R(1)=0, ~~R(g)=-\xi R(e)-\eta R(f),~~\beta_2\neq 0,~~\xi\neq 0~\text{ or }~ \nu\neq 0.$$
\end{lemma}

\begin{proof}
	Let  $\delta_0\neq 0$. Then $v,u$ is a basis  of $\ker R$. We have
	$$0=[v,u]=[\delta_0 g,\beta_0 e+\gamma_0 f]=2\delta_0\beta_0 e-2\delta_0\gamma_0 f. $$
	Therefore,  $u=\alpha_0 1$. Then we may suppose that $u=1$ and
	$v= g+\xi e+\eta f$. In this case  ${\rm Im} R=span(R(e), R(f))$ and
	$ R(e), R(f)$ is  a basis of  ${\rm Im} R$.
	
	Assume  $R(g)=0$. Then $v=g$.  Since by (\ref{kercon1}) $[R(e)+\lambda e,g]\in \ker R\cap (\Bbbk e + \Bbbk f)=0$  we have  $[R(e)+\lambda e, g]=0$. Similarly,
	$[R(f)+\lambda f,g]=0$.
	
	Let $ R(e)=\alpha 1+\sigma g+\beta e+\gamma f$. Then $ -2\beta
	e+2\gamma f-2\lambda e=0$ and  we get $\beta=-\lambda$,
	$\gamma=0$. Consequently,  $$ R(e)=\alpha 1+\sigma g-\lambda e.$$
	Similarly we obtain  $ R(f)=\alpha_1 1+\sigma_1 g-\lambda f.$ Then
	$$2(\lambda\sigma f+\lambda\sigma_1 e)=-[\sigma g,\lambda f]-[\lambda e,\sigma_1
	g]=$$ $$[R(e),R(f)]=R([R(e),f]+[e,R(f)])= -2R(\sigma
	f+\sigma_1e)=$$ $$-2\sigma (\alpha_1 1+\sigma_1
	g-\lambda f)- 2\sigma
	_1(\alpha 1+\sigma  g-\lambda e).$$ From these equalities we obtain that
	$\sigma \alpha_1=-\sigma _1\alpha$ and $\sigma \sigma _1=0$.
	
	Moreover, if ${\rm Im} R$ is an abelean algebra, then  $\sigma
	=\sigma _1=0$. In this case  $$R(1)=0, ~~R(g)=0, ~~R(e)=\alpha 1-
	\lambda e,~~R(f)=\alpha_1 1-\lambda f.$$
	
	If ${\rm Im} R$ is a non-abelean algebra, then  $\sigma\neq 0$ or $\sigma_1\not=0$. In this case
	$$R(1)=0, ~~R(g)=0, ~~R(e)=\alpha 1+ \sigma g-\lambda e,~~R(f)=-\lambda f,$$
	or
	$$R(1)=0, ~~R(g)=0, ~~R(e)= -\lambda e~~,R(f)=\alpha_1 1+\sigma g -\lambda f.$$
	These cases coincide up to isomorphism  $\phi$.
	
	Assume  $g\not \in \ker R$. Then $(\Bbbk \cdot e+\Bbbk \cdot f)\cap \ker R=0$
	and $\xi\not= 0$ or $\eta\not= 0$. Therefore,
	$$H_4= (\Bbbk \cdot e+\Bbbk \cdot f) \oplus \ker R.$$ Let
	$$
	a=\alpha e+\beta f+\sigma 1+\gamma v
	$$ and $[a,v]=0$. Then
	$$0=[a,v]=[\alpha e+\beta f,v]=[\alpha e+\beta f,g]=-2
	\alpha e+2\beta f.$$ Consequently, $a\in \ker R$.
	
	For any  $a\in H_4$,  by (\ref{kercon1}), we have $$[R(a)+\lambda a, v]\subseteq (\Bbbk \cdot e+\Bbbk \cdot f  )\cap\ker R=0.$$  From this equality we get
	$R(a)+\lambda a\in \ker R$. Hence $R(R(a))=-\lambda R(a)$ for any
	$a\in H_4$ and $R(a)=-\lambda a $ for any $a\in {\rm Im }R$.

	Let $R(e)=-\lambda e+\alpha_11+\beta_1 v$, $R(f)=-\lambda
	f+\alpha_2 1+\beta_2 v$. Then
	$$[R(e),R(f)]=R([R(e),f]+[e,R(f) ])=-2\beta_1 R(f)-2\beta_2 R(e).$$
	
	If ${\rm Im}R$  is an abelean algebra, then $\beta_1=\beta_2=0.$
	Consequently,
	$$R(1)=0, ~~R(e)=-\lambda e+\alpha 1,~~R(f)=-\lambda f+\alpha_1 1,$$
	$$ R(g)=\xi\lambda e+\eta \lambda f-(\xi\alpha+\eta\alpha_1)1,~~ \xi\neq 0~\text{ or }~ \nu\neq 0.$$
	If ${\rm Im}R$ is non-abelean, then $\beta_1\neq 0$, or $\beta_2\neq 0$
	and $\beta_1\alpha_2=-\beta_2\alpha_1$, $\beta_1\beta_2=0.$
	
	Thus,  $$
	R(e)=-\lambda e,~~R(f)=-\lambda f+\alpha_2 1+\beta_2 (g+\xi e+\eta f),
	$$
	$$R(1)=0,~~ R(g)=-\xi R(e)-\eta R(f),~~\beta_2\neq 0,\xi\neq 0~\text{ or }~ \nu\neq 0,
	$$  or
	$$  R(e)=-\lambda e+\alpha_1e +\beta_1 (g+\xi e+\eta f),  ~~R(f)=-\lambda f,$$
	$$R(1)=0, ~~R(g)=-\xi R(e)-\eta R(f),~~\beta_1\neq 0,~~\xi\neq 0~\text{ or } ~\nu\neq 0.$$
	These cases coincide up to isomorphism  $\phi$.
\end{proof}

\begin{lemma} If $\delta_0=0$, then the operator  $R$ has  one from the  following forms.

If  $f\not\in \ker R$, then

{\rm (i)} $$R(1)=-\gamma 1, ~~R(g)=-\lambda g+w,~~ w\in I,  ~~R(e)=-\eta 1,~~R(f)=1, ~~\eta\neq 0.$$

{\rm (ii)}  $$ R(f)=\lambda e,~~R(g)=-\lambda g+w_1,~~w_1\in I,~~R(e)=-\lambda e,~~R(1)=-\lambda e.$$

{\rm (iii)}  $$R(f)=-\lambda f,~~R(g)=-\lambda g+w_1~~,w_1\in I,~~R(e)=\lambda f,~~R(1)=\gamma \lambda f.$$

{\rm (iv)} $$R(1)=-\gamma R(f),~~R(g)=\zeta R(f) +\mu 1=-\zeta\lambda f+(\zeta\xi+\mu) 1 +\zeta\nu e,$$
$$R(e)=0, ~~R(f)=-\lambda f+\xi 1 +\nu e,~~\mu \neq 0.$$

{\rm (v)}  $$ R(f)= g+\beta_1e-(\lambda+\sigma_1)f+\xi 1,R(1)=-\gamma R(f),$$ 	$$R(e)=0, R(g)=\sigma_1 R(f)+\mu 1,\, \mu\neq 0.$$

{\rm (vi)} $$R(1)=-\gamma 1, ~~R(g)=\sigma_1g+\sigma_1(\beta_1e+\gamma_1f)+\xi 1, ~~   R(e)=0,R(f)=1.$$

{\rm (vii)} $$  R(1)=0,~~R(g)=\xi \lambda f+\alpha_21 +\beta_2e,$$
$$R(e)=0, ~~R(f)=-\lambda f+\alpha_11+\beta_1e,~~\beta_2\neq \xi\beta_1.$$

{\rm (viii)} $$ R(1)=\lambda \gamma f, ~~R(g)=-\lambda g+w_1,~~w_1\in I, ~~R(e)=0, ~~R(f)=-\lambda f.$$

{\rm (ix)} $$R(1)=0, ~~ R(g)=\sigma_1g  +w_1, w_1\in I,~~ R(e)=0,~~R(f)=-\lambda f.$$

{\rm (x)} $$R(1)=0,  ~~R(g)=-\frac{\lambda}{2}g+w_1, ~~w_1\in I,~~R(e)=0, ~~R(f)=-\lambda^{-1} e.$$

If $f\in \ker R$, then

{\rm (xi)} $$ R(1)=\sigma g+\alpha 1+w, ~~w\in \ker B,  ~~R(g)=\xi R(1)+\mu 1,$$
$$ R(e)=R(f)=0 , ~~ \sigma\neq 0,\,\mu\neq 0. $$

{\rm (xii)} $$R(1)=\xi R(g)+\mu 1,~~ R(g)=\sigma g+\alpha 1+w, ~~w\in \ker B, $$
$$  R(e)=R(f)=0 ,\, \sigma\neq 0,~~\mu\neq 0 .$$

{\rm (xiii)}  $$R(1)=\alpha_1 1+\beta_1 e+\gamma_2f,~~R(g)=\alpha_2 1+\beta_2e +\gamma_2f,~~R(e)=R(f)=0.$$
\end{lemma}

\begin{proof} Since $\delta_0=0$,  we may suppose that $u=\alpha_0 1+\beta_0 e+\gamma_0 f$,
	$v = e+\eta f$ is a basis of $\ker R$.
	
	Let us consider the case  $R(f)\neq 0$. In this case  ${\rm Im }R=span(R(f), R(g))$,   $R(f), R(g)$ is a basis of ${\rm Im} R$ and
	$\alpha_0\neq 0$.  Therefore, we may suppose that $u=1+\gamma f.$
	
	Since $H_4=  \Bbbk g+I $, then $$R(f)=\sigma g+w,~~R(g)=\sigma_1 g+w_1,$$
	where  $w\in I$ and $w_1\in I$.
	We have $R(e)=-\eta R(f)$. Hence, by (\ref{YBop}),
	$$0=[R(e), R(f)]=R(-\eta\sigma[g,f]+\sigma[e,g])=2R(\sigma\eta f-\sigma e).$$
	Therefore, $\sigma\eta R(f)=0$. Consequently,
	$\sigma\eta=0$.
	
	Similarly $$[R(e),R(g)]=R(-\eta[w,g]+\sigma_1[e,g]+\lambda [e,g])=$$
	$$-\eta[R(f),R(g)]=-\eta R([w,g]+\sigma_1[f,g]+\lambda [f,g]).
	$$
	From this equality  we get  $(\sigma_1+\lambda)\eta R(f)=0$.  Then
	$(\sigma_1+\lambda)\eta=0$.

	Consider the first case  $\eta\neq 0$. In this case $\sigma=0$ and  $\sigma_1=-\lambda$.
	
	We have  ${\rm Im}R\cap I\neq 0$. Then  $0\neq a=\alpha_1 1+\beta_1
	e+\gamma_1 f\in {\rm Im}R$. Since  ${\rm Im}R$ is a subalgebra then
	$$4\lambda^2( \beta_1 e+\gamma_1f)=[R(g), ~~[R(g),a]]\in {\rm Im}R.$$
	  Therefore,  $\alpha_1 1, \beta_1 e+\gamma_1f\in {\rm Im}R$. Consequently,  $\alpha_1=0$,  or    $\beta_1 =\gamma_1=0$ otherwise
	the elements $\alpha_1 1,$ $\beta_1 e+\gamma_1f$, $-\lambda g+w_1$ are
	linearly dependent.
	
	Let  $1\in {\rm Im}R$. Then $1,w$ are  linearly dependent and we may suppose that $R(f)=1$.
	
	Thus,
	$$R(f)=1, ~~R(g)=-\lambda g+w,~~ w\in I,~~ R(1)=-\gamma 1, ~~R(e)=-\eta 1,~~ \eta\neq 0.$$
	
	Let  $0\neq\beta_1 e+\gamma_1f\in{\rm Im}R$. Then the algebra  ${\rm
		Im}R$ is non-abelean and
	$$2\lambda(\beta_1 e-\gamma_1f)=[R(g),\beta_1 e+\gamma_1f]\in{\rm
		Im}R.$$ From here we get $\beta_1e,\gamma_1f\in {\rm
		Im}R$. Consequently, $\beta_1\neq 0$ or  $\gamma_1\neq 0. $
	
	Let $\beta_1\neq 0$.  Then $$R(f)=\alpha e, R(g)=-\lambda g+w_1$$
	for  some  $\alpha\in \Bbbk$ and  $\alpha\neq 0$. Moreover, by
	(\ref{YBop}), we have
	$$[R(f),R(g)]=R(-2\alpha e -2\lambda f+2\lambda
	f)=2\alpha^2\eta e.$$ On the other hand,
	$[R(f),R(g)]=2\alpha\lambda e.$   Consequently,
	$\alpha=\frac{\lambda}{\eta}$.
	Therefore, $$R(f)=\frac{\lambda}{\eta}e,~~R(g)=-\lambda g+w_1,~~w_1\in I,~~R(e)=-\lambda e,~~R(1)=-\frac{\gamma\lambda}{\eta}e, ~~\eta\neq 0.$$
 Thus, we may assume  $$R(f)=\lambda e,~~R(g)=-\lambda g+w_1,~~w_1\in I,~~R(e)=-\lambda e,~~R(1)=-\lambda e.$$

	Let $\gamma_1\neq 0$.  Then
	$$R(f)=\alpha f, ~~R(g)=-\lambda g+w_1$$
	for some $\alpha\in \Bbbk$ and $\alpha\neq 0$. Moreover, by
	(\ref{YBop}), we have
	$$[R(f), R(g)]=R(2\alpha f -2\lambda f+2\lambda f)=2\alpha^2 f.$$
	On the other hand,
	$[R(f),R(g)]=-2\alpha\lambda f.$   Consequently,
	$\alpha=-\lambda$.
	Thus, we may assume $$R(f)=-\lambda f,~~R(g)=-\lambda g+w_1~~,w_1\in I,~~R(e)=\lambda f,~~R(1)=\gamma \lambda f.$$

	Now consider the second case $\eta=0.$ Then  $u= 1+\gamma f$, $v=e.$
	Since ${\rm Im} R \cap J\neq 0$ then $0\neq b =\alpha_1 g +\beta_1 e+\gamma_1 f \in {\rm Im} R$.
	
	Let $1\in {\rm Im} R$. If ${\rm Im} R$ is non-abelean, then
	$[R(f), R(g)]\neq 0$. Therefore,
	$${\rm Im} R\cap (\Bbbk e + \Bbbk f)\neq 0.$$ Then
	${\rm Im} R\subseteq I$ and ${\rm Im} R$ is abelean.
	
	If $\sigma=\sigma_1=0$ then, by  (\ref{YBop}),
	$$0=[R(g),R(f)]=R([w_1,f]+[g,w]-2\lambda f)=R([g,w]-2\lambda f).$$
	From here we get  $w=-\lambda f+\xi 1 +\nu e$, $\xi,\nu \in \Bbbk$ since  $R(f)\neq 0$.
	Thus,
	$$R(1)=-\gamma R(f),~~R(e)=0, ~~R(f)=-\lambda f+\xi 1 +\nu e,  $$
	$$ R(g)=\zeta R(f) +\mu 1=-\zeta\lambda f+(\zeta\xi+\mu) 1 +\zeta\nu e, ~~ \mu \neq 0.$$

	If $\sigma\neq 0$ then   $\alpha_1\neq 0$ since  ${\rm Im}R$ is abelean. Therefore, we may suppose   $b=g +\beta_1 e+\gamma_1 f.$
	
	Since $[R(f),b]=[R(g),b]=0$ then $w=\sigma(\beta_1 e+\gamma_1
	f)+\xi 1$. Moreover, $R(g)=\frac{\sigma_1}{\sigma}R(f)+\mu 1$ for some non-zero $\mu\in \Bbbk$.
	Consequently,  $$ R(f)=\sigma g+\sigma(\beta_1 e+\gamma_1 f)+\xi
	1,$$ $$ R(g)=\frac{\sigma_1}{\sigma}R(f)+\mu 1,\mu\neq 0.$$
	
	By (\ref{YBop}) we get
	$$0=[R(g),R(f)]=R([R(g),f]+[g,R(f)]-2\lambda f)=$$
	$$R(-2\sigma_1 f+2\sigma\beta_1 e -2\sigma\gamma_1f-2\lambda f).$$
	From here $ \gamma_1=-\frac{\lambda+\sigma_1}{\sigma}$ and
	$$ R(f)=\sigma g+\sigma\beta_1e-(\lambda+\sigma_1)f+\xi 1,R(1)=-\gamma R(f),$$
	$$R(e)=0, R(g)=\frac{\sigma_1}{\sigma}R(f)+\mu 1,\,\sigma \neq 0, \mu\neq 0.$$
 Thus, we may assume
	$$ R(f)= g+\beta_1e-(\lambda+\sigma_1)f+\xi 1,R(1)=-\gamma R(f),$$
	$$R(e)=0, R(g)=\sigma_1 R(f)+\mu 1,\, \mu\neq 0.$$
	
	If $\sigma=0$, $\sigma_1\neq 0$ then $R(f)=\frac{\sigma}{\sigma_1}R(g)+\mu1=\mu1$, where $\mu\neq 0$. We may suppose that $R(f)=1$.
	Thus  $$R(g)=\sigma_1g+\sigma_1(\beta_1e+\gamma_1f)+\xi 1, R(f
	)=1, R(e)=0, B (1)=-\gamma1.$$

	Let  $1\not \in{\rm Im}R$.   Suppose     ${\rm Im}R$ is an abelean  algebra. Let us show that in this case  $\sigma=\sigma_1=0$. Indeed, suppose that, for example, $\sigma\neq 0$. There are exist
	$\xi,\zeta,\mu \in \Bbbk$ such that $0\neq \mu 1+ y=\xi R(f)+\zeta
	R(g)$, where $y \in \Bbbk e+ \Bbbk  f$, since $\ker B\cap I\neq 0$. Then we have
	$$4\sigma y=[R(f),[R(f),\mu 1+ y]]=0.$$  Consequently, $y=0$ and
	$0\ne \mu 1\in {\rm Im}B$.
  Therefore,  $\sigma=\sigma_1=0$ and
	$R(f)=w$, $R(g)=w_1$.

 By  (\ref{kercon}), $[R(g)+\lambda g,u ]\in
	\ker B$.
	Since $[R(g)+\lambda g,u ]=-2\lambda\gamma f$, we have   $1\in \ker B$.
	
	Let  $w=\alpha_11+\beta_1e+\gamma_1f.$ Then, by   (\ref{YBop}),
	$$0=[R(g),R(f)]=R([g,w]-2\lambda f).$$
	From here we get  $\gamma_1=-\lambda $. Therefore,
	$$R(f)=\alpha_11+\beta_1e-\lambda f.$$ Note that
	$\ker R\cap {\rm Im}R\neq 0$. Hence,   $$R(g)=\alpha_21
	+\beta_2e-\xi\lambda f.$$ The elements $1, f,e$ are linear independent.
	Since  $1\not\in {\rm Im}B$ then  $\beta_2\neq \xi\beta_1$ or
	$\xi\alpha_1=\alpha_2$.
	Thus,
	$$ R(f)=-\lambda f+\alpha_11+\beta_1e, R(1)=0,\,
	R(e)=0,$$
	$$R(g)=\alpha_21 +\beta_2e-\xi \lambda f, \beta_2\neq \xi\beta_1 \text{ or }\xi\alpha_1=\alpha_2.$$

	Let ${\rm Im} R$ be non-abelean. Then    $[R(g),R(f)]\neq  0$ and
	$\sigma \neq 0$ or $\sigma_1 \neq 0$ and $1\not\in {\rm Im}R$.
	Therefore, $[R(g),R(f)]=\xi f$ or
	$[R(g),R(f)]=\xi e$,  where $\xi\neq 0$. In the last case we can assume $[R(g),R(f)]= e$.
	
	Let $w=\alpha_11+\beta_1 e+\gamma_1f$.  From  (\ref{YBop}) we have
	$$
	[R(g),R(f)]=R(-2\sigma_1f+2\beta_1e-2\gamma_1f-2\lambda f ) = -2(\sigma_1+\gamma_1+\lambda )R(f).
	$$
	Consequently,
	$\sigma_1+\gamma_1+\lambda\neq 0$, $\sigma=0$, $\alpha_1=0$ and
	$\sigma_1\neq 0$.
	
	If $[R(g),R(f)]=\xi f$,  then $\beta_1=0$ and  $R(f)=\gamma_1f$, $\gamma_1\neq 0$.
	In this case we get   $\sigma_1\gamma_1=(\sigma_1+\gamma_1+\lambda)\gamma_1$.
	From here  $\gamma_1=-\lambda$. It means $R(f)=-\lambda f$.
	Since $R(1)=-\gamma R(f)$ and $[R(1),R(g)]=B[R(1), g]$ then  $(\sigma_1+\lambda)\gamma=0$.
	Thus,
	$$R(f)=-\lambda f, ~~R(g)=-\lambda g+w_1, ~~w_1\in I,~~ R(1)=\lambda\gamma f, ~~R(e)=0,
	$$
	or
	$$
	R(f)=-\lambda f, ~~R(g)=\sigma_1g  +w_1,~~w_1\in I, ~~R(1)=0, ~~R(e)=0,
	$$

	If $[R(g),R(f)]=e$, then  $\gamma_1=0$, $R(f)=\beta_1e$ and $\beta_1\neq 0$, since  $R(f)\neq 0$.
	Then
	$$
	2\sigma_1\beta_1e=[R(g),~~R(f)]=-2(\sigma_1+\lambda )\beta_1e.
	$$
	From here $\sigma_1=-\frac{\lambda}{2}.$ Since  $[R(g),R(f)]=e$ we have $\beta_1=-\lambda^ {-1}.$
	By (\ref{kercon}),  $[R(g)+\lambda g, u]\in \ker R.$
	Since  $[R(g)+\lambda g, u]=-\gamma \lambda f$ and  $f\not\in\ker R$ we have $1\in \ker   R$.
	Thus,
	$$
	R(f)=-\lambda^ {-1}e, ~~R(g)=-\frac{\lambda}{2}g+w_1,~~ w_1\in I,~~R(1)=R(e)=0.
	$$

	Let $R(f)=0 $. Then $R(e)=0$.  In this case
	$$[R(H_4),R(H_4)]\subseteq R([H_4,H_4])\subseteq R(\Bbbk e+\Bbbk f)=0.$$
	Therefore, the algebra  ${\rm Im} R$ is abelean. Moreover, ${\rm
		Im} R=span(R(1), R(g))$ and $R(1), R(g)$ is a basis of ${\rm Im} R$.
	
	Let
	$$
	R(1)=\sigma_1 g_1+w_1, ~~R(g)=\sigma_2 g+w_2,~~\sigma_1,~~\sigma_2\in \Bbbk, ~~ w_1,~~ w_2\in I.
	$$
	Assume that $\sigma_1\neq  0$ or  $\sigma_2 \neq  0$. In this case  ${\rm Im} R\cap K=0$, else $[{\rm
		Im}R, {\rm Im}R]\neq 0$  and $1\in {\rm Im} R\cap I$. Then
	$$
	R(e )=R(f)=0, ~~R(1)=\sigma g+\alpha 1+w, ~~w\in \ker B,~~ \sigma\neq 0,
	$$
	$$  R(g)=\xi R(1)+\mu 1,~~\mu\neq 0, $$ or
	$$
	R(e )=R(f)=0, ~~R(g)=\sigma g+\alpha 1+w, ~~w\in \ker B,~~ \sigma\neq 0,
	$$
	$$  R(1)=\xi R(g)+\mu 1,~~ \mu\neq 0 .$$
	
	Assume $\sigma_1=\sigma_2=0$. Then
	$$R(1)=\alpha_1 1+\beta_1 e+\gamma_2f, ~~ R(g)=\alpha_2 1+\beta_2e +\gamma_2f,~~ R(e)=R(f)=0.$$
\end{proof}

\subsection{The kernel $\ker R$ is nonabelean.} Let $\ker R$ be an  nonabelean algebra.
Since $\ker R\cap I\neq  0$ then $0\neq v=\alpha1+\beta e+\gamma
f\in \ker R$.

Let  $0\neq u=\delta  1+\xi g+\mu e+\nu f\in \ker R.$ We can suppose that $\xi=1$.  Moreover,
$$2\beta e-2\gamma f=[u,v]\in \ker B, 4\beta e+4\gamma f=[u, [u,v]]\in \ker R .$$
From here we get  $\alpha1,\beta e, \gamma f\in \ker R$. The
subalgebra $\ker R$ is nonabelean, hence $\alpha=0$ and we may
suppose that  $e\in \ker R$, $f\not\in \ker R$ and $u=\delta
1+g+\gamma f$. Consequently, $e$ and $u$ is a basis of  $\ker R$. In
this case $R(1),R(f)$ is a basis of  ${\rm Im}R$.

 The elements $1,f,u,e$ are a basis of  $H_4$.
  Let $R(f)=\sigma_1
1+\sigma_2f+\sigma_3 u+\sigma_4e$ and  $h=\delta 1+g$. Since
$R(h)=-\gamma R(f)$ then, by (\ref{YBop}),
$$
0=[R(h), R(f)]=R([R(h), f]+[h, R(f)]+\lambda[h, f])=-2R(\sigma_2 f+\lambda f).
$$
Consequently,  $\sigma_2=-\lambda$
and  $$R(f)=\sigma_1 1-\lambda f+\sigma_3 u+\sigma_4e.$$

Let $R(1)=\delta_11+\delta_2f+\delta_3u+\delta_4e.$ By
(\ref{YBop}), $$[R(1), R(f)]=R([R(1),f])=-2\delta_3R(f).$$ Similarly
$$2\delta_3\gamma
R(f)=2\delta_3[R(1),R(h)]=R([R(1),h])=2R(\delta_2f+\gamma\delta_3f-\delta_4e).$$
From here $\delta_2=0$ and $$R(1)=\delta_11+\delta_3u+\delta_4e.$$

Since  ${\rm Im} R\cap I\neq  0$ then  $0\neq w=\alpha 1+\beta e+\gamma f\in {\rm Im}R$.

Let  $\sigma_3\neq 0$. We can suppose that $\sigma_3=1$. Then
$$[R( f),w]=2\beta e-2\gamma
f\in {\rm Im} R.$$ From here  $\alpha 1,\beta e, \gamma f\in
{\rm Im} R$. The elements $R(f),1,e,f$ are linearly independent.
Therefore,  ${\rm Im} R\neq span(1,e,f)$.

Let $1\in  {\rm Im} R$. Then ${\rm Im} R =span(R(f),1)$, ${\rm Im}
R$ is abelean and  $-\lambda f+ u+\sigma_4e\in {\rm Im} R$,
$\delta_3u+\delta_4e\in {\rm Im} R$. Also $$2\lambda\delta_3
f+2(\sigma_4\delta_3-\delta_4)e=[R(1),-\lambda f+ u+\sigma_4e]\in
{\rm Im} R.$$ From here $\delta_3=0$, $\delta_4=0$, else  ${\rm Im}R=span (1,e, f)$. Consequently, $R(1)=\delta_11$ and    $\delta_1\neq 0$.

Thus, $$R(1)=\delta_11, R(f)=\sigma_11+g+(\gamma-\lambda) f+\sigma_4e,$$
$$ R(g)=-\delta R(1)-\gamma R(f), R(e)=0, \delta_1\neq 0.$$

If $e\in {\rm Im}R$ then $R(1)-\delta_4e=\delta_11+\delta_3u\in {\rm
	Im }R$ and
$$2\lambda\delta_3 f+2\delta_3\sigma_4e=[\delta_11+\delta_3u,R(f)]\in {\rm Im }R.$$
From here $\delta_3=0$.  Then we may suppose  $\delta_1=0$,
$R(1)=\delta_4 e$  and $\delta_4\neq 0$.

Thus,  $$R(1)=\delta_4e, R(f)=\sigma_11+g +(\gamma-\lambda) f+\sigma_4e,$$ $$
R(g)=\delta R(1)-\gamma R(f), R(e)=0,\delta_4\neq 0.$$

If  $f\in {\rm Im} R$ then $$2(\delta_4-\sigma_4\delta_3)e=[R(f)+\lambda f, R(1)]\in {\rm Im}R.$$
From here  $\delta_4-\sigma_4\delta_3=0$. Then
$ \delta_3u+\delta_4e=\delta_3(u+\sigma_4e)$ and $\delta_3\neq 0$, else
$R(1)=\delta_1 1\in {\rm Im} R$. Therefore,
$$(\delta_1-\delta_3\sigma_1)1=R(1)-\delta_3 (R(f) +\lambda f)\in {\rm Im}R.$$
From here $\delta-\delta_3\sigma_1=0$ and $R(1)=\delta_3(R(f)+\lambda f).$

Thus,
$$R(1)=\delta_3(R(f)+\lambda f), R(f)=\sigma_11 +g+(\gamma-\lambda) f+\sigma_4e,$$
$$ R(g)=\delta R(1) -\gamma R(f), R(e)=0,\delta_3 \neq 0.$$

Let $\sigma_3=0$ and $\delta_3\neq 0$. Then  $R(f)=\sigma_11-\lambda
f+\sigma_4 e$ and
$$[R(1),R(f)]=2\lambda\delta_3f+2\delta_3 \sigma_4e.$$ Therefore, $f\in
{\rm Im}R$, $\sigma_4e\in {\rm Im}R$, $\sigma_11\in {\rm Im}R$.
From here $\sigma_1=\sigma_4=0$.

Thus,  $$R(1)=\delta_1 1+\delta_3u+\delta_4 e, R(f)=-\lambda
f, R(g)=-\delta R(1)-\gamma R(f), R(e)=0.$$

Let $\sigma_3=\delta_3=0$. Then $$R(1)=\delta_11+\delta_4e,
R(f)=\sigma_11-\lambda f+\sigma_4e,$$ $$ R(g)=-\delta R(1) -\gamma R(f), R(e)=0, \delta_1\neq 0 \text{ or  }
\delta_4\neq 0 .$$

Holds

\begin{lemma}
	Let $\dim\ker R=2$ and let $\ker R$ be an nonabelean algebra. Then
	the operator  $R$ is one of following
	
	\noindent {\rm (i) } $$R(1)=\delta_11, R(f)=\sigma_11+g+(\gamma-\lambda)
	f+\sigma_4e,$$ $$ R(g)=-\delta R(1)-\gamma R(f), R(e)=0,
	\delta_1\neq 0.$$ {\rm (ii) } $$R(1)=\delta_4e, R(f)=\sigma_11+g
	+(\gamma-\lambda) f+\sigma_4e,$$ $$ R(g)=-\delta R(1)-\gamma R(f),
	R(e)=0,\delta_4\neq 0.$$ {\rm (iii)} $$R(1)=\xi(R(f)+\lambda f),
	R(f)=\sigma_11 +g+(\gamma-\lambda) f+\sigma_4e,$$
	$$ R(g)=-\delta R(1)-\gamma R(f), R(e)=0,\xi\neq 0.$$
	{\rm (iv)} $$R(1)=\delta_1 1+\delta_3g+\gamma\delta_3 f+\delta_4 e,
	R(f)=-\lambda f,$$
	$$ R(g)=-\delta R(1)-\gamma R(f), R(e)=0, \delta_3\neq 0.$$ {\rm
		(v)}
	$$R(1)=\delta_11+\delta_4e, R(f)=\sigma_11-\lambda f+\sigma_4e,$$ $$
	R(g)=-\delta R(1)-\gamma R(f), R(e)=0, \delta_1\neq 0 \text{ or }
	\delta_4\neq 0 .$$
\end{lemma}

\bigskip


\section{The kernel has the dimension one}

Let $\dim\ker R=1$. Then  $\dim {\rm Im}R=3$ and ${\rm Im}R\cap
K\neq 0$. Moreover,
$$\dim({\rm Im}R\cap I)\geq 2, ~~\dim({\rm Im}R\cap
J)\geq 2.
$$

Assume ${\rm Im}R$ is an  abelean algebra. Then  $[{\rm Im}R\cap
K,{\rm Im}R\cap J]= 0$. Therefore, $\Bbbk e+\Bbbk f \subseteq {\rm Im}R $,
since $\dim({\rm Im}R\cap J)\geq 2$. From here we get  ${\rm
Im}R=I$.

Let $R(1)=\alpha 1+\beta e+\gamma f$, $R(e)=\alpha_1 1+\beta_1 e+\gamma_1f$,
$R(f)=\alpha_2 1+\beta_2 e+\gamma_2 f$. By (\ref{YBop}),
$$0=[R(g),R(e)]=2R((\beta_1+\lambda )e-\gamma_1f).$$ From here
$(\beta_1+\lambda )e-\gamma_1f\in \ker B$. Similarly
$\beta_2e-(\gamma_2+\lambda )f\in \ker B$ and $\beta e-\gamma f\in
\ker B.$

If ${\rm Im}R\cap\ker R=0$, then
$$R(1)=\alpha 1, R(e)=\alpha_11-\lambda e,R(f)=\alpha_21-\lambda f, \alpha\neq 0,$$
$$R(g)=\beta R(1)+\delta R(e)+\gamma R(f) .$$

If ${\rm Im}R\cap\ker R\neq 0$, then $ 0\neq w=\xi 1+\eta e+\mu f\in \ker R.$

Assume  $\xi=0$. Then we may assume  that  $w=e+\eta f\in \ker R$.
Consequently, $$R(1)=\alpha 1+\beta e-\beta\eta f,R(e)=\alpha_1 1+\beta_1e-(\beta_1+\lambda)\eta f,$$ $$
R(f)=\alpha_2 1+\beta_2 e-(\beta_2\eta +\lambda)f, R(e)=-\eta R(f).$$
Then  $\eta=0,\alpha_1=0,\beta_1=0$.
Thus,  $$R(1)=\alpha 1+\beta e , R(e)=0, R(f)=\alpha_2 1+\beta_2e-\lambda f,$$
$$ R(g)=\sigma  1+\delta e+\nu f, \beta(\nu\alpha_2+\sigma\lambda)-\alpha(\delta\lambda+\nu\beta_2)\neq 0
.$$

Let $\xi\neq 0$. Since, by  (\ref{kercon}), $[R(g)+\lambda g,w]\in
\ker B$, then $2\eta e-2\mu f\in\ker B$. Therefore, we may suppose that
$w=1$.
Thus,  $$R(1)=0,R(e)=\alpha 1-\lambda e, R(f)=\alpha_1
1-\lambda f,$$ $$R(g)=\sigma  1+\mu e+\nu f,\lambda\sigma
+\alpha\mu+\alpha_1\nu\neq  0.$$

Hence, we have

\begin{theorem} \label{th0} Let $1,e,f$ be a basis of  ${\rm Im}R$. Then each Rota--Baxter operator $R$ acts on the basis   $1$, $g$, $e$, $f$ by one of the following rules:

{\rm (i)}
$$R(1)=\alpha 1, R(e)=\alpha_11-\lambda e,R(f)=\alpha_21-\lambda f, \alpha\neq 0,$$
$$R(g)=\beta R(1)+\delta R(e)+\gamma R(f) .$$
{\rm (ii)}
$$R(1)=\alpha 1+\beta e , R(e)=0, R(f)=\alpha_2 1+\beta_2e-\lambda f,$$
$$ R(g)=\sigma  1+\delta e+\nu f, \beta(\nu\alpha_2+\sigma\lambda)-\alpha(\delta\lambda+\nu\beta_2)\neq 0
.$$
{\rm (iii)} $$R(1)=0,R(e)=\alpha 1-\lambda e, R(f)=\alpha_1
1-\lambda f,$$ $$R(g)=\sigma  1+\mu e+\nu f,\lambda\sigma
+\alpha\mu+\alpha_1\nu\neq  0.$$
\end{theorem}

Assume that the algebra  ${\rm Im} R$ is not abelean.    By Lemma \ref{lm1}, the vector space ${\rm Im}R$ has
one of the  following bases
 $\alpha
1+g,ef$ or $1,g+\gamma f,e$ or $1,g+\gamma e,f.$

\begin{theorem} \label{th1}  Let $h=\alpha 1+g,e,f$ be a basis of  ${\rm Im}R$. Then each Rota--Baxter operator  $R$ acts on the basis $1,h,e,f$ by one of the following rules:

    \noindent {\rm (i)}     $$R(1)=\xi h+\mu e+\nu f,R(h)=\xi^{-1}\nu\beta_3e+\xi^{-1}\nu\lambda f,$$ $$ R(e)=0, R(f)=\beta_3e-\lambda f, \xi\nu\beta_3\neq 0.$$
{\rm (ii)}$$R(1)=\mu e+\nu f, R(h)=-\lambda h+\beta_1e+\gamma_1f,$$ $$R(e)=-\lambda e, R(f)=-\lambda e,\nu\neq 0.$$
{\rm (iii)} $$R(1)=0, R(h)=\alpha_1h+\beta_1e+\gamma_1f,R(e)=-\lambda e,$$ $$R(f)=-\lambda f,\alpha_1\neq 0,-\lambda. $$ {\rm (iv)} $$R(1)=\mu e+\nu f, R(h)=-\lambda h +\beta_1e+\gamma_1f,R(e)=-\lambda e,R(f)=-\lambda f. $$ {\rm (v)}
$$R(1)=\xi h+\mu e+\nu f,R(e)=-\lambda e,R(f)=-\lambda f,$$ $$ R(h)=\alpha_1h +\xi^{-1}(\alpha_1+\lambda)\mu e+\xi^{-1}(\alpha_1+\lambda)\nu f. $$ {\rm (vi)}  $$R(1)=0, R(h)=-\frac{\lambda}{2}h+\beta_1 e+\gamma_1f,$$ $$R(e)=-\lambda e,R(f)=\beta_3e-\lambda f,\beta_3\neq 0.$$  {\rm (vii)}
$$R(1)=\xi h+\mu e+\nu f, R(h)=-\frac{\lambda}{2}h+\xi^{-1}(\frac{\lambda}{2}\mu+\beta_3\nu) e+\xi^{-1}\frac{\lambda}{2}\nu f,$$ $$R(e)=-\lambda e,R(f)=\beta_3e-\lambda f,\beta_3\neq 0.$$
{\rm  (viii)} $$ R(1)=0,  R(h)=-\frac{\lambda}{2}h+\beta_1e+\gamma_1f,$$
$$ R(e)=\beta_2e+f, R(f)=\beta_2(\beta_2+\lambda)e+\beta_2f,\beta_2\neq 0,-\lambda.$$
\end{theorem}

\begin{proof}
    The algebra  ${\rm Im}R$ has   the  same multiplication table for the basis  $h=\alpha 1+g,e,f$ as in Lemma  \ref{lm2}.

    Let
$$R(1)=\xi h+\mu e+\nu f, R(h)=\alpha_1h+\beta_1e+\gamma_1 f, $$
$$R(e)=\alpha_2h+\beta_2e+\gamma_2 f,R(f)=\alpha_3h+\beta_3e+\gamma_3 f $$

The case  1.1  of Lemma \ref{lm2} impossible since $\ker R> 1$.

In the case  1.2 $\alpha_2\neq0$ from $[R(1),R(e)]=R([R(1),e])$ we get
$$\xi\beta_2 e-\xi\gamma_2 f-\alpha_2\mu e+\alpha_2\nu f=\xi R(e).$$
From here $\xi =0, \mu=0, \nu=0$, i.e. $R(1)=0$.
Then $\dim\ker R>1$.

In the case  1.4 from
$[R(h), R(1)]=R([h,R(1)])$ and $[R(f), R(1)]=$ $R([f,R(1)])$  we get
$$\lambda\mu+2\xi\beta_1=2\nu\beta_3, 2\xi\gamma_1=-\lambda \nu,\xi \beta_3=0.$$
From here    $\xi=0$ since  $\beta_3\neq 0$. Then $\mu=\nu=0$. Therefore  $R(1)=0$ and
$\dim\ker  R>2$,  since $R(e)=0$.

In the case 1.5 from $[R(h), R(1)]=R([h, R(1)])$ we get  $$
\xi\beta_1=\alpha_1\mu, \xi\gamma_1=(\lambda+\alpha_1)\nu.$$
If  $\xi=0$, then $\mu\neq 0$, else  $\dim\ker R>1$.
Then  $\alpha_1=0$,   $\nu =0$ and  $\dim\ker R>1$.  It means  $\xi\neq 0$. Then
$\dim\ker R>1$.

In the case  1.6  from $[R(h), R(1)]=R([h, R(1)])$ we get  $$
\xi\beta_1=\nu\beta_3, \xi\gamma_1=\lambda\nu.$$ If  $\xi=0$ then $\nu=0$ and    $\mu\neq0$, else  $\dim\ker R>1$. Then   $\dim\ker R>1$. Therefore,  $\xi \neq 0$. From here  $\nu\neq 0$, else  $\dim\ker R>1$.
Thus, $$R(1)=\xi h+\mu e+\nu f,R(h)=\xi^{-1}\nu\beta_3e+\xi^{-1}\nu\lambda f,$$ $$ R(e)=0, R(f)=\beta_3e-\lambda f, \xi\nu\beta_3\neq 0.$$

In the case 1.7 we have $\beta_3\neq0$. From $[R(h),R(1)] =R([h,R(1)])$ and  $[R(h),R(1)] =R([h,R(1)])$ we get $$\xi\beta_1=0,\, \xi\gamma_1=(\lambda+\beta_3)\nu,\, \xi\beta_3=0.$$ From here $\xi= 0$. Then  $\nu\neq 0$, else  $\dim\ker R>1$. Consequently,  $\beta_3=-\lambda.$

Thus,  $$R(1)=\mu e+\nu f, R(h)=-\lambda h+\beta_1e+\gamma_1f,$$ $$R(e)=-\lambda e, R(f)=-\lambda e,\nu\neq 0.$$

In the case  1.8 from  $[R(h),R(1)] =R([h,R(1)])$ we get
$$\xi\beta_1=(\alpha_1+\lambda)\mu, \xi\gamma_1=(\alpha_1+\lambda)\nu.$$

Thus, $$R(1)=0, R(h)=\alpha_1h+\beta_1e+\gamma_1f,R(e)=-\lambda e,$$
$$R(f)=-\lambda f,\alpha_1\neq 0,-\lambda $$ or
$$R(1)=\mu e+\nu f, R(h)=-\lambda h +\beta_1e+\gamma_1f,$$ $$R(e)=-\lambda e,R(f)=-\lambda f $$
or  $$R(1)=\xi h+\mu e+\nu f,R(e)=-\lambda e,R(f)=-\lambda f,$$
$$ R(h)=\alpha_1h +\xi^{-1}(\alpha_1+\lambda)\mu e+\xi^{-1}(\alpha_1+\lambda)\nu f. $$

In the case  1.9 from  $[R(h),R(1)] =R([h,R(1)])$ we get
$$2\xi\beta_1=\lambda\mu+2\beta_3\nu,2\xi\gamma_1=\lambda\nu.$$

Thus,  $$R(1)=0, R(h)=-\frac{\lambda}{2}h+\beta_1 e+\gamma_1f,$$
$$R(e)=-\lambda e,R(f)=\beta_3e-\lambda f,\beta_3\neq 0$$ or
$$R(1)=\xi h+\mu e+\nu f, R(h)=-\frac{\lambda}{2}h+\xi^{-1}\left(\frac{\lambda}{2}\mu+\beta_3\nu\right) e+\xi^{-1}\frac{\lambda}{2}\nu f,$$
$$R(e)=-\lambda e,R(f)=\beta_3e-\lambda f,\beta_3\neq 0.$$

In the case 1.10 $\alpha_2=\alpha_3=0, \beta_3=\gamma_2^{-1}\beta_2(\beta_2+\lambda)\neq 0$.
From $[R(e), R(1)]=R([e,R(1)])$ we get  $\xi=0$.
From  $[R(h), R(1)]=R([h,R(1)])$ we get
$$\mu\gamma_2=\nu \left(\beta_2+\frac{\lambda}{2}\right), \nu\beta_3=\mu\left(\beta_2+\frac{\lambda}{2}\right).$$
If $\mu\nu\neq 0$ then  $$\beta_2(\beta_2+\lambda)=\left(\beta_2+\frac{\lambda}{2}\right)^2.$$ It is contradiction. It  means $\mu\nu=0$. Then  $\mu=\nu=0$.
Therefore,  $$R(1)=0,  R(h)=-\frac{\lambda}{2}h+\beta_1e+\gamma_1f,$$
$$R(e)=\beta_2e+\gamma_2f, R(f)=\gamma_2^{-1}\beta_2(\beta_2+\lambda)e+\beta_2f,\beta_2\neq 0,-\lambda,\gamma_2\neq 0 .$$

Thus, $$R(1)=0,  R(h)=-\frac{\lambda}{2}h+\beta_1e+\gamma_1f,$$
$$R(e)=\beta_2e+f, R(f)=\beta_2(\beta_2+\lambda)e+\beta_2f,\beta_2\neq 0,-\lambda .$$
\end{proof}

\begin{theorem} \label{th2} Let $ 1, h=g+\gamma f,e$ be a basis of  ${\rm Im}R$.
Then each Rota--Baxter operator  $R$ acts on  the basis  $1, h, e, f$ by one of the following rules:

\noindent {\rm (i)} $$R(1)=\gamma_31+\beta_3e,
R(h)=\gamma_11-\frac{\lambda}{2}h+\beta_1e,$$ $$ R(e)=0, R(f)=\mu
e,\mu \neq 0,\gamma_3\neq 0.$$
{\rm (ii)}
$$R(1)=\gamma_3 1, R(h)=\gamma_1
1+\alpha_1h+\beta_1e,$$ $$
R(e) =-\lambda e, R(f)=0, \alpha_1\neq 0,-\lambda, \gamma_3\neq 0.$$
{\rm (iii)}
$$R(1)=\gamma_3 1, R(h)=\gamma_11-\lambda
h+\beta_1e,$$ $$ R(e) =-\lambda e, R(f)=\nu 1+\mu e, \gamma_3\neq 0
\text{ or  } \nu\neq 0.$$   {\rm (iv)}
$$R(1)=\gamma_3 1, R(h)=\gamma_11-\lambda h,R(e) =-\lambda e,$$ $$
R(f)=\nu1+\xi h+\mu e,\gamma_3\neq 0 \text{ or } \nu\neq 0 \text{ or
} \xi\neq 0.$$
{\rm (v)}
$$R(1)=\gamma_31+\beta_3e,    R(h)=\gamma_1 1-\lambda h+\beta_1e,R(e)=-\lambda e,$$
$$  R(f)=\nu 1+\mu f,\beta_3\neq 0, \gamma_3 \neq 0 \text{ or } \nu\neq 0.$$
{\rm (vi)}$$R(1)=\gamma_31+\alpha_3h+\beta_3e,
R(h)=\gamma_11+\alpha_1h+\beta_3\alpha_3^{-1}(\alpha_1+\lambda)e,$$ $$R(e)=-\lambda e, R(f)=0,\alpha_1\neq 0, \gamma_3\alpha_1-\gamma _1\alpha_3\ne0.$$
{\rm (vii)}$$R(1)=\gamma_31+\alpha_3h+\lambda^{-1}\alpha_3e,
R(h)=\gamma_11+\beta_1e,R(e)=\gamma_21-\lambda e,$$
$$ R(f)=0, \alpha_3(\beta_1\gamma_2+\lambda\gamma_1)\neq 0.$$
{\rm (viii)}$$R(1)=\gamma_31, R(h)=\gamma_11-\lambda h-\xi^{ -1}\lambda\mu e, R(e)=0,$$
$$R(f)=\nu1+\xi h+\mu e, \gamma_3\neq 0,\xi\neq0, \gamma_1\xi+\lambda\nu\neq 0. $$ {\rm (ix)} $$ R(1)=\gamma_31, R(h)=\gamma_11-\lambda h, R(e)=\gamma_2 1,$$
$$R(f)=\nu1+\xi  h, \gamma_2\neq 0,\xi\neq0. $$
\end{theorem}

\begin{proof} The basis  $h, e,1$  of the algebra  ${\rm Im} R$ has the same
multiplication table as in Lemma~\ref{lm3}.

Let $$R(h)=\alpha_1 h+\beta_1e+\gamma_1
1,\, R(f)=\xi h+\mu e+\nu 1.$$  From (\ref{YBop}) for the element   $h, f $ we get
\begin{equation}\label{h-gx-1}(\alpha_1\mu-\xi\beta_1) e=\mu R(e)-(\alpha_1+\lambda) R(f).
\end{equation}

Let  $$R(1)=\alpha_3 h+\beta_3e+\gamma_3 1.$$  Then from
(\ref{YBop})  for the element   $1, f$  we get
\begin{equation}\label{z-gx-1}(\alpha_3\mu-\xi\beta_3)e=-\alpha_3R(f).
\end{equation}

In the case   2.1 of  Lemma \ref{lm3} $R(e) =0$ and $\alpha_1\neq
-\lambda$. Then from (\ref{h-gx-1}) we get
$$\xi=0, \nu=0,(2\alpha_1+\lambda)\mu=0.$$
Since
$\dim\ker R=1$ then  $R(f)\neq 0$. Therefore,   $\mu\neq 0$ and
$\alpha_1=-\frac{\lambda}{2}$.  From  (\ref{z-gx-1})
we get  $\alpha_3=0$.
Thus, $$R(1)=\gamma_31+\beta_3e,
R(h)=\gamma_11-\frac{\lambda}{2}h+\beta_1e,$$ $$ R(e)=0, R(f)=\mu
e,\mu \neq 0,\gamma_3\neq 0.$$

In the cases  2.2 and 2.3   of Lemma \ref{lm3} we have $\alpha_1=0$
and $R(e)=0$. Therefore, by  (\ref{h-gx-1}),
$$\xi =0, \nu=0,\mu=0.$$ Then $R(f)=0 $ and $\dim \ker R>1$.

In the case  2.4.  $\alpha_3=\beta_3=0$. From  (\ref{h-gx-1}) we get
$$(\alpha_1+\lambda)\xi= (\alpha_1+\lambda)\nu=0,2(\alpha_1+\lambda)\mu=\xi\beta_1.$$

If $\alpha_1\neq -\lambda$ then $\xi=\nu=0$, $\mu =0$.

Thus, $$R(1)=\gamma_ 1, R(h)=\gamma_1
1+\alpha_1h+\beta_1e,$$ $$
R(e) =-\lambda e, R(f)=0, \alpha_1\neq 0,-\lambda, \gamma_3\neq 0.$$

If $\alpha_1=-\lambda$, then  $\xi\beta_1=0$.

Thus,   $$R(1)=\gamma_3 1, R(h)=\gamma_11-\lambda h+\beta_1e,$$ $$
R(e) =-\lambda e, R(f)=\nu 1+\mu e, \gamma_3\neq 0 \text{ or  }
\nu\neq 0$$ or
$$R(1)=\gamma_3 1, R(h)=\gamma_11-\lambda h,R(e) =-\lambda e,$$ $$
R(f)=\nu1+\xi h+\mu e,\gamma_3\neq 0 \text{ or } \nu\neq 0 \text{ or } \xi\neq 0.$$

In the case 2.5  $\alpha_1=-\lambda, \alpha_3=0$. From (\ref{h-gx-1}) and (\ref{z-gx-1}) we get    $\xi\beta_1=0$ and
$\xi\beta_3=0$.

By the previous case, we may suppose $\xi=0$, $\beta_3\neq 0$.

Thus,
$$R(1)=\gamma_31+\beta_3e,    R(h)=\gamma_1 1-\lambda h+\beta_1e,R(e)=-\lambda e,$$
$$  R(f)=\nu 1+\mu f,\beta_3\neq 0 \text { and } \gamma_3 \neq 0 \text{ or } \nu\neq 0.$$

In the case  2.6 $\alpha_3\neq 0$.   From    (\ref{z-gx-1}) we get
$$\alpha_3\xi=\alpha_3\nu=0,   2\alpha_3\mu =\xi\beta_3.$$
From here  $\xi=\mu=\nu=0$.

Thus,  $$R(1)=\gamma_31+\alpha_3h+\beta_3e,
R(h)=\gamma_11+\alpha_1h+\beta_3\alpha_3^{-1}(\alpha_1+\lambda)e,$$
$$R(e)=-\lambda e, R(f)=0,\alpha_1\neq 0, \gamma_3\alpha_1-\gamma _1\alpha_3\ne0.$$

In the case 2.7 $\alpha_1=0.$  From (\ref{h-gx-1}) we get
$$\xi\beta_1=2\mu\lambda, \gamma_2\mu=\lambda\nu, \xi=0.$$ From here  $\mu=\nu=0$.
Since  $\dim\ker R=1$ then $(\gamma_1\lambda+\gamma_2\beta_1)\alpha_3\neq 0$.

Thus,  $$R(1)=\gamma_31+\alpha_3h+\lambda^{-1}\alpha_3e,
R(h)=\gamma_11+\beta_1e,R(e)=\gamma_21-\lambda e,$$
$$ R(f)=0, (\gamma_1\lambda+\gamma_2\beta_1)\alpha_3\neq 0.$$

In the case 2.8  $\alpha_1=-\lambda$ and $\beta_3=-\lambda^
{-1}\alpha_3\beta_1$.  From (\ref{h-gx-1}) and (\ref{z-gx-1} ) we
get
$$\lambda\mu=-\xi\beta_1, \mu\gamma_2=0,\alpha_3\xi=\alpha_3\nu=0,2\alpha_3\mu=\xi\beta_3.$$

Let  $\alpha_3\neq 0$. Then  $\xi=\mu=\nu=0.$ From here   $\dim\ker
R>1$. Therefore,  $\alpha_3=0$ and $\beta_3=0$.

Let   $\xi= 0$. Then $\mu=0$ and  $\dim\ker R>1$. Therefore,  $\xi\neq
0$.

Thus  $$R(1)=\gamma_31, R(h)=\gamma_11-\lambda h-\xi^{ -1}\lambda\mu e, R(e)=0,$$
$$R(f)=\nu1+\xi h+\mu e, \gamma_3\neq 0,\xi\neq0, \gamma_1\xi+\lambda\nu\neq 0$$ or
$$ R(1)=\gamma_31, R(h)=\gamma_11-\lambda h, R(e)=\gamma_2 1,$$
$$R(f)=\nu1+\xi  h, \gamma_2\neq 0,\xi\neq0. $$

In the  case 2.9 we consider the basis
$h_1=h-\alpha_2^{-1}\alpha_1e, e,1$ of   ${\rm Im}R$. Then
$$R(h_1)=(\gamma_1-\alpha_2^{-1}\alpha_1\gamma_2)1,
R(1)=\gamma_31,$$ $$ R(e)=\alpha_2h-\lambda e+\gamma_21.$$

Let $R(f)=\xi h_1+\mu e +\nu 1.$ Then from (\ref{YBop})  for the elements
$h_1,f$ we get $\lambda R(f)=\mu R(e).$ From here  $dim\ker R>1$.
\end{proof}

\section{The kernel has the dimension zero}

In this section we describe nondegenerate Rota--Baxter operators of non-zero weight  on the Lie algebra $H_4^{(-)}$
\begin{theorem}
    Let $R$ be a nondegenerate Rota--Baxter operator of non-zero weight on $H_4^{(-)}$.

    If   $R(I)\neq I$ then for some $\alpha\in \Bbbk$ the operator  $R$ has on the the basis  $1,h=\alpha 1+g,e,f$  the following form
     $$R(1)=\xi  h+\mu e+\nu f,R(h)=\alpha_11+\beta_1h+\xi^{-1}(\lambda+\beta_1)\mu e+\xi^{-1}(\lambda+\beta_1)\nu
       f,$$ $$ R(e)=-\lambda e, B(f)=-\lambda f,\xi\alpha_1\neq 0.$$

    If $R(I)=I$, then for some $\alpha\in \Bbbk$ and   $h=\alpha 1+g$
    the operator $R$ has one of the following forms

    \noindent {\rm (i)} $$R(1)=\xi 1+\mu e +\nu f , R(g)=-\lambda h+\beta_1e+\gamma_1 f,$$ $$
    R(e)=-\lambda e, B(f)=\lambda f,\xi\neq 0.$$
{\rm (ii)} $$R(1)=\xi 1 , B(g)=\alpha_1 h+\beta_1e+\gamma_1 f,$$   $$
R(e)=-\lambda e, R(f)=-\lambda f,\alpha_1\neq 0, -\lambda, \xi\neq 0.$$
    {\rm (iii)} $$R(1)=\xi 1, R(g)=-\frac{\lambda}{2}h+\beta_1e+\gamma_1 f,
    R(e)=-\lambda e, R(f)=\beta_3e
    -\lambda  f, \beta_3\xi\neq 0 .$$ {\rm (iv)} $$ R(1)=\xi 1 ,R(g)=-\frac{\lambda }{2}h+\beta_1e+\gamma_1f,R(e)=\beta_2e +f,$$
 $$ R(f)=\beta_2(\beta_2+\lambda)e+\beta_2 f, \xi\neq 0, \beta_2\neq 0, -\lambda, $$
\end{theorem}

\begin{proof}   Let $R(I )\neq I$. Then by Lemma  \ref{lm1}, $R(I)=J_\alpha$ for some  $\alpha\in \Bbbk$.

    For any  $u,v\in I$ we have  $$[R(u),R(v)]=R([R(u),v]+[u,R(v)]).$$ Therefore,  $R^{-1} \colon J_\alpha \to I$ is derivation.

    Let  $h=\alpha1+g$. Since  $
[g, e] =[h,e]=2e$ and  $e\in J_\alpha\cap I$   then
       $$2R^{-1}(e) =R^{-1}([h,e])=[R^{-1}(h),e]+[h,R^{-1}(e)]=[g,R^{-1}(e)].$$ Consequently,
     $R(e)=\mu e$, $\mu\in \Bbbk$. Similarly  $R(f)=\nu f$, $\nu\in \Bbbk$.

   Let   $R(g)=\alpha_11+\beta_1g+\gamma_1e+\delta_1 f$. Then
    $$2\beta_1e=[R(g),e]=\mu^{-1}[R(g),R(e)]=\mu^{-1}R([R(g),e]+[g,R(e)]+2\lambda e)=$$
    $$2\beta_1e+2\mu e+2\lambda e.$$  From here we get  $\mu=-\lambda$. Similarly $\nu=-\lambda.$

  Put   $R(1)=\xi h+\mu e +\nu f$ and $R(h)=\alpha_11+\beta_1h+\gamma_1e+\delta_1 f.$
    Then $\xi\neq 0$, else $1+\lambda^{ -1}\mu e+\lambda^{ -1}\nu f\in ker R$.
    By (\ref{YBop}), for $1,h$  we get
    $$\gamma_1 =\xi^{-1}(\lambda+\beta_1)\mu, \delta_1=\xi^{-1}(\lambda+\beta_1)\nu.$$

    Thus,  $$R(1)=\xi  h+\mu e+\nu f,R(h)=\alpha_11+\beta_1h+\xi^{-1}(\lambda+\beta_1)\mu e+\xi^{-1}(\lambda+\beta_1)\nu  f,$$
    $$ R(e)=-\lambda e, R(f)=-\lambda f. $$
    Since $R$ is a nondegenerate operator then  $\xi\alpha_1\neq 0$.

    Let  $R(I)=I$. Then there are  non-zero   $\alpha\in \Bbbk$ and $w\in I$ such that $g=R (\alpha g+w)$. Let $e=R(u)$, $u\in I$. By (\ref{YBop}),
    $$2e=[g,e]=[R(\alpha g+w),R(u)]=$$ $$=R([R(\alpha g+w),u]+[\alpha g+w, R(u)]+\lambda[\alpha g+w,u]).$$
    From here  $e=\sigma R(e)+\delta R(f)$ since  $[H_4,H_4]\subseteq \Bbbk e+\Bbbk f.$

    Similarly  $f=\xi R(e)+\mu R(f).$ Then $R(\Bbbk e+\Bbbk f)\subseteq \Bbbk e+\Bbbk f$ since $R$ is a nondegenerate operator.

    The vector space   $R(J)$ is an ideal of  $H_4^{(-)}$.
    Indeed, $$[R(J),H_4^{(-)}]\subseteq [R(J),R(H_4^{(-)})]\subseteq R(\Bbbk e+\Bbbk f)\subseteq R(J).$$
    Therefore, by Lemma \ref{lm1}, $R(J)=J_\alpha$ for some  $\alpha\in \Bbbk$.
    Consequently,  $R(g)=\alpha_1(\alpha 1+g)+\beta_1e+\gamma_1f$ and $\alpha_1\neq 0.$

    Let  $\overline{R} \colon  J\to J$ be the mapping, defined by the rule
    $$\overline{R}(g)=R(g)-\alpha_1\alpha1,      \overline{R}(e)=R(e),\overline{R}(f)=R(f).$$
   The mapping $\overline{R}$ is nondegenerate on  $J$.  Since $1$ lies in the center of the  algebra  $H_4^{(-)}$ and $[H_4,H_4]\subseteq \Bbbk e+\Bbbk f$ then  $\overline{R}$ is a Rota-Baxter operator on $J$.
    Then, by Lemma  \ref{lm2}, only the cases  1.8, 1.9 and 1.10 hold for $\overline{R}$.

    Put $R(1)=\xi 1+\mu e+\nu  f$ and  $h=\alpha 1+g$.

In the cases 1.8 and 1.9
    $$ R(g)=\alpha_1 h+\beta_1e+\gamma_1f,R(e)=-\lambda e, R(f)=\beta_3e-\lambda f,$$ $\alpha_1\neq 0, \beta_3=0$
     or  $\alpha_1=-\frac{\lambda}{2},\beta_3\neq 0.$ Then, by  (\ref{YBop}), we get
    $$-2\alpha_1\mu e+2\alpha_1\nu f=[R(1),R(g)]=R([R(1),g])=R(-2\mu e+2\nu f).$$
From here  $(\alpha_1+\lambda)\mu=\beta_3\nu$, $(\alpha_1+\lambda)\nu=0$.

    In the case  1.8 we get  $$R(1)=\xi 1+\mu e +\nu f , R(g)=-\lambda h+\beta_1e+\gamma_1 f,$$ $$
    R(e)=-\lambda e, B(f)=-\lambda f,\xi\neq 0,$$ or
$$R(1)=\xi 1 , R(g)=\alpha_1 h+\beta_1e+\gamma_1 f,$$   $$
R(e)=-\lambda e, R(f)=-\lambda f,\alpha_1\neq 0, -\lambda, \xi\neq 0.$$
    In the case  1.9 we get   $$R(1)=\xi 1, R(g)=-\frac{\lambda}{2}h+\beta_1e+\gamma_1 f,
    R(e)=-\lambda e, R(f)=\beta_3 e-\lambda  f,\beta_3 \xi\neq 0 .$$

 In the case  1.10 $$ R(1)=\xi 1 ,R(g)=-\frac{\lambda }{2}h+\beta_1e+\gamma_1f,R(e)=\beta_2e +\gamma_2f,$$
 $$R(f)=\gamma_2^{-1}\beta_2(\beta_2+\lambda)e+\beta_2 f, \xi\neq 0, \beta_2\neq 0, -\lambda, \gamma_2 \neq 0, $$

 Thus, we may assume
 $$ R(1)=\xi 1 ,R(g)=-\frac{\lambda }{2}h+\beta_1e+\gamma_1f,R(e)=\beta_2e +f,$$
 $$R(f)=\beta_2(\beta_2+\lambda)e+\beta_2 f, \xi\neq 0, \beta_2\neq 0, -\lambda, $$
 \end{proof}

\begin{question}
We have found Lie algebra RB-operator on $H^{(-)}_4$ and in the paper \cite{Ma} were  found associative algebra RB-operator on the associative algebra $(H_4, \cdot, 1)$.
What Lie algebra RB-operators on $H^{(-)}_4$  are not associative RB-operators on $(H_4, \cdot, 1)$?
\end{question}

\bigskip


\section*{Acknowledgments}
Authors are grateful to  V.~Gubarev, M.~Goncharov, A. Pozhidaev,  and P. Kolesnikov for the fruitful discussions and useful suggestions.
Authors are also grateful to participants of the seminar ``\'{E}variste Galois''
at Novosibirsk State University for attention to our work.

Valery G. Bardakov is supported by the Russian Science Foundation (RSF
 24-21-00102) for work in sections 2 and 5.

Viktor N. Zhelyabin is supported by the Program of Fundamental Research RAS,
	(FWNF-2022-0002) for work in sections 3 and 4.


\medskip

\noindent Valeriy G. Bardakov \\
Sobolev Institute of Mathematics,  Acad. Koptyug ave. 4, 630090 Novosibirsk, Russia; \\
Novosibirsk State Agrarian University,
Dobrolyubova str., 160, 630039 Novosibirsk; \\
Regional Scientific and Educational Mathematical Center of Tomsk State University; \\
Lenin ave. 36, 634009 Tomsk, Russia \\
email: bardakov@math.nsc.ru

\medskip
\noindent  Viktor N. Zhelaybin \\
Sobolev Institute of Mathematics,  Acad. Koptyug ave. 4, 630090 Novosibirsk, Russia; \\
email: vicnic@mail.math.nsc.ru

\medskip
\noindent Igor Nikonov \\
Lomonosov Moskow  State University \\
email: nikonov@mech.math.msu.su


\begin{thebibliography}{67}

\bibitem{AnBai} H. An, C. Bai, \textit{From Rota-Baxter Algebras to Pre-Lie Algebras}, J. Phys. A. 1 (2008),
015201, 19 p

\bibitem{Atkinson}
F. V. Atkinson, \textit{Some aspects of Baxter's functional equation}, J. Math. Anal. Appl., 7 (1963) 1--30.


\bibitem{BN}
V. G. Bardakov, I. M. Nikonov, \textit{Relative Rota--Baxter operators on groups and Hopf algebras}, arXiv:2311.09311v1.


\bibitem{BG-2} V.~G.~Bardakov, V. Gubarev, \textit{Rota--Baxter operators on groups}, Proc. Indian Acad. Sci. (Math. Sci.), 133, no. 4 (2023).

\bibitem{BG-1} V.~G.~Bardakov, V. Gubarev, \textit{Rota--Baxter groups, skew left braces, and the Yang--Baxter equation},  J. Algebra,  596 (2022), 328--351.


\bibitem{Baxter}
G. Baxter, \textit{An analytic problem whose solution follows from a simple algebraic identity}, Pacific J. Math.,  10 (1960) 731--742.

\bibitem{BNZH} V. G. Bardakov, I. M. Nikonov, and V. N. Zhelyabin  \textit{ Hopf algebra and group Rota--Baxter operators,} in progress.

\bibitem{BGP} P. Benito, V. Gubarev, A. Pozhidaev, \textit{Rota--Baxter operators on quadratic algebras,}
Mediterr. J. Math. 15 (2018), 23 p. (N189).

\bibitem{BG} T. A. Bolotina, V. Gubarev, \textit{Rota--Baxter Operators on the Simple Jordan Superalgebra $D_t$},
Siberian Mathematical Journal, (63) 4 (2022)  637-650.




\bibitem{Cartier} P. Cartier, \textit{On the structure of free Baxter algebras}, Adv. Math., 9 (1972) 253--265.

\bibitem{Drinfeld} V. G. Drinfel'd,
\textit{ Hamiltonian structures on Lie groups, Lie bialgebras and the geometric meaning of classical Yang-Baxter equations}, Dokl. Akad. Nauk SSSR (268) 2
(1983) 285--287.

\bibitem{DuBaiGou} C. Du, C. Bai, L. Guo, \textit{3-Lie bialgebras and 3-Lie classical Yang-Baxter equations
in low dimensions}, Linear Multilinear A. 66 (8) (2018) 1633--1658.


\bibitem{GonGub} M. Goncharov, V. Gubarev,
\textit{Rota--Baxter operators of nonzero weight on the matrix algebra of order three}, Linear and Multilinear Algebra,  (7) 6 (2022) 1055-1080

\bibitem{Goncharov} M.E. Goncharov, \textit{
The description of Rota-Baxter operators of nonzero weight on general linear complex Lie algebra of order 2}, Siberian Electronic Mathematical Reports,   (19) 2 (2022) 870-879.
\bibitem{GK} M. E. Goncharov,  D. E. Kozhukhar’, \textit{ Rota--Baxter Operators of Nonzero Weight on a Complete Linear Lie Algebra of Order Two},
Algebra and Logic, (61) 1  (2022) 67-70.
\bibitem{Goncharov2021} M. Goncharov, \textit{Rota–Baxter operators on cocommutative Hopf algebras}, J. Algebra, 582, no. 1 (2021),
39--56.
\bibitem{GuoMonograph}
L. Guo, \textit{An Introduction to Rota--Baxter Algebra}, Surveys of Modern Mathematics, vol. 4,
International Press, Somerville (MA, USA); Higher education press, Beijing, 2012.
\bibitem{GuoLShen} L. Guo, H. Lang, and Yu. Sheng, \textit{Integration and geometrization of Rota--Baxter Lie algebras}, Adv.
Math., 387 (2021).

\bibitem{Kolesniukov} P.S. Kolesnikov, \textit{Homogeneous averaging operators on simple finite conformal Lie
algebras}, J. Math. Phys. 56 (2015) 071702, 10 p.
24
\bibitem{Konovalova} E. I. Konovalova, \textit{Double Lie algebras}, Ph.D. Thesis, Samara State University, 2009.
189 p. (in Russian).


\bibitem{Ma}
 T. S. Ma, A. Makhlouf, S. Silvestrov, \textit{Rota-Baxter cosystems and coquasitriangular mixed bialgebras},J. Algebra Appl.,  20 (2021), 2150064.


\bibitem{Ma-1}
Tianshui Ma, Jie Li, Liangyun Chen, Shuanhong Wang, \textit{Rota--Baxter operators on Turaev's Hopf group (co)algebras I:
Basic definitions and related algebraic structures}, J. Geometry and Physics,  175 (2022) 104469.

\bibitem{Pan} Yu Pan, Q. Liu, C. Bai, L. Guo, \textit{PostLie algebra structures on the Lie algebra $sl(2, C)$}, Electron. J. Linear Algebra 23 (2012) 180--197.

\bibitem{PeiBai} J. Pei, C. Bai, and L. Guo, \textit{Rota-Baxter operators on sl(2, C) and solutions of the
classical Yang--Baxter equation}, J. Math. Phys. 55 (2014), 021701, 17 p.

\bibitem{Rota} G.-C. Rota, \textit{Baxter algebras and combinatorial identities}, I, Bull. Amer. Math. Soc.
75 (1969) 325--329.

\bibitem{STyanSh} M. A. Semenov-Tyan-Shanskii, \textit{What is a classical r-matrix?}, Funct. Anal. Appl. 17
(1983) 259–272.

\bibitem{Sokolov} V. V. Sokolov, \textit{Classification of constant solutions of the associative Yang--Baxter
equation on $Mat_3$}, Theor. Math. Phys. (3) 176 (2013) 1156--1162.

\bibitem{Tricomi}
F.~G.~Tricomi, \textit{On the finite Hilbert transformation}, Quart. J. Math., no. 2 (1951),  199--211.

\bibitem{T}
 V. G. Turaev, \textit{Quantum Invariants of Knots and 3-Manifolds}, De Gruyter Studies in Mathematics, 2016, 596 pp.

\bibitem{ZLMZ}
Huihui Zheng, Fangshu Li,  Tianshui Ma, Liangyun Zhang,  \textit{Hopf brace, braid equation and bicrossed coproduct}, arXiv:1912.01392v1.


\end{thebibliography}
\end{document}